\documentclass[12pt,a4paper]{article} 
\usepackage{amsmath,pstricks,latexsym,theorem,epsfig} 
\usepackage{amsfonts,amssymb,latexsym,epsfig,amsmath}
 
\usepackage{amsmath,latexsym,theorem} 
\newtheorem{Theorem}{Theorem}[section] 
\newtheorem{Definition}{Definition}[section] 
\newtheorem{Proposition}{Proposition}[section] 
\newtheorem{Lemma}{Lemma}[section] 
\newtheorem{Corollary}{Corollary}[section] 
 
\theorembodyfont{\rmfamily} 
\newtheorem{Remark}{Remark}[section]

\newcommand{\rec}[1]{{(\ref{#1})}} 
\def \R{I\!\!R}

\def \Z{Z\!\!\!Z} 
 
\def\11{1\!\!1}

\def\rec#1{{(\ref{#1})}}

\newcommand{\ba}{\begin{array}} 
\newcommand{\ea}{\end{array}}

\begin{document} 
 
 \title{\bf Compactness and Bubbles Analysis for $1/2$-harmonic Maps}
\author{ Francesca Da Lio}
\date{} 

\maketitle 
\begin{abstract} In this paper we study compactness and quantization properties of sequences of   $1/2$-harmonic maps $u_k\colon\R\to {\cal{S}}^{m-1}$ such that $\|u_k\|_{\dot H^{1/2}(\R,{\cal{S}}^{m-1})}\le C\,.$ More precisely we show that there exist a weak $1/2$-harmonic map
$u_\infty\colon\R\to {\cal{S}}^{m-1}$, a possible empty set  $\{a_1,\ldots,a_\ell\}$ in $\R$ such that 
up to subsequences
$$(|(-\Delta)^{1/4}u_k|^2 \rightharpoonup |(-\Delta)^{1/4}u_{\infty}|^2)dx+\sum_{i=1}^{\ell}\lambda_i \delta_{a_i},~~\mbox{in Radon measure}\,,$$
as $k\to +\infty$,  with $\lambda_i\ge 0\,.$

The convergence of $u_k$ to $u_\infty$  is strong in $\dot W^{1/2,p}_{loc}(\R\setminus\{a_1,\ldots,a_\ell\})$, for every $p\ge 1\,.$  We     quantify  the loss of energy in the weak convergence and   we show that in the case of  non-constant $1/2$-harmonic maps with values in   $   {\cal{S}}^2\,$ one has $\lambda_i=2 \pi n_i$, with $n_i$ a positive  integer\,.

 \end{abstract}
 {\small {\bf Key words.} Fractional harmonic maps, nonlinear elliptic PDE's, regularity of solutions, commutator estimates.}\par
 {\small { \bf  MSC 2000.}  58E20, 35J20, 35B65, 35J60, 35S99}
 
\tableofcontents
\section{Introduction}
In the paper \cite
{DLR1} Rivi\`ere and the author started the investigation of the following 
$1$-dimensional quadratic Lagrangian
\begin{equation}\label{nonloclagr}
L(u)=\int_{\R}|(-\Delta)^{1/4} u(x)|^2 dx\,,
\end{equation}
where $u\colon\R\to{\cal{N}}$, ${\cal{N}}$ is  a smooth   $k$-dimensional submanifold of  $\R^m$   which is at least $C^2$, compact and without boundary.
We observe that \rec{nonloclagr} is a simple  model of Lagrangian which is  invariant under the trace of conformal maps that keep invariant the half space $\R^2_+$: the M\"obius group.
\par
  Precisely let  $\phi\colon\R^2_+\to \R^2_+$ in $W^{1,2}(\R^2_+,\R^m)$ be a conformal map of degree $1$, i.e. it satisfies
  \begin{equation}\label{conf}
\left\{\begin{array}{l}
\displaystyle{|\frac{\partial\phi}{\partial x}|=|\frac{\partial\phi}{\partial y}|}
\\
\displaystyle{\langle \frac{\partial\phi}{\partial x},\frac{\partial\phi}{\partial y}\rangle=0}\\
{\rm det}\nabla\phi\ge 0 ~~{\mbox and} ~~\nabla\phi\ne 0\,.
 \end{array} \right.
\end{equation}
Here   $\langle\cdot,\cdot\rangle$ denotes the
standard Euclidean inner product in $\R^m$\,.\par

  We denote by
  $\tilde \phi$ the restriction of $\phi$ to  $\R$. Then we have
  $L(u\circ \tilde\phi)=L(u)\,.$
  
Moreover  $L(u)$ in \rec{nonloclagr}  coincides with the semi-norm $\|u\|^2_{\dot{H}^{1/2}(\R)}$ and  the following identity holds
 \begin{eqnarray}\label{char}
\int_{\R}|(-\Delta)^{1/4} u(x)|^2 dx=\inf\left\{\int_{\R^2_+}|\nabla \tilde u|^2 dx:~\tilde u\in W^{1,2}(\R^2,\R^m),~~\mbox{trace $\tilde u=u$}\right\}\,.
\end{eqnarray}

The Lagragian $L$ extends to map $u$ in the following function space
 $$
 \dot{H}^{1/2}(\R,{\cal{N}})=\{u\in \dot{H}^{1/2}(\R,\R^m):~~u(x)\in {\cal{N}}, {\rm a.e},\}\,.
 $$
The operator $(-\Delta)^{1/4} $ on ${\R}$ is defined by means of the the Fourier transform as follows  $$\widehat{(-\Delta)^{1/4} u}=|\xi|^{1/2}\hat u\,,$$
 (given a function f, $\hat f$ denotes the Fourier transform of $f$).\par

We denote by  $\pi_{\cal{N}}$ the orthogonal projection from $\R^m$  onto ${\cal{N}}$ which happens to be a $C^\ell$ map in a sufficiently small
neighborhood of ${\cal{N}}$ if ${\cal{N}}$ is assumed to be $C^{\ell+1}$. We now introduce the notion of $1/2$-harmonic map into a manifold.

 \begin{Definition}
 \label{df1}
 A map $u\in \dot{H}^{1/2}(\R,{\cal{N}})$   is called a weak $1/2$-harmonic map into $\cal{N}$ if for any $\phi\in  \dot{H}^{1/2}(\R,\R^m)\cap L^{\infty}({\R},{\R}^m)$ there holds
 \[
 \frac{d}{dt}L(\pi_{\cal{N}}(u+t\phi))_{|_{t=0}}=0\quad.
 \]
 \end{Definition}
 In short we say that a weak $1/2-$harmonic map is a {\it critical point of $L$ in  $\dot{H}^{1/2}(\R,{\cal{N}})$ for perturbations in the target}.

\medskip
We next give some geometric motivations related to the study of the problem \rec{nonloclagr}\,.\par
First of all  variational problems of the form \rec{nonloclagr} appear as a simplified model of 
  renormalization area in hyperbolic spaces, see for instance \cite{AM}\,. 
  There are also some   geometric connections which are being investigated  in the paper \cite{DLR3}  between $1/2$-harmonic maps and the so-called free boundary sub-manifolds and optimization problems of eigenvalues. 
 With this regards we refer the reader also to the papers \cite{FS1,FS2}\,.
  Finally  $1/2-$harmonic maps into the circle ${\cal{S}}^1$   might appear for instance in the asymptotic of equations in phase-field theory for fractional reaction-diffusion such as
 \[
\epsilon^2\,(-\Delta)^{1/2}u+u (1-|u|^2)=0\,,
 \]
 where $u$ is a complex valued ''wave function''.\par

\medskip

In this paper we consider the case ${\cal{N}}={\cal{S}}^{m-1}$. It can be shown (see, \cite{DLR1}) that every weak  $1/2$-harmonic map satisfies the following Euler-Lagrange equation
\begin{equation}\label{ELE}
(-\Delta)^{1/2}u\wedge u=0 ~~\mbox{in ${\cal{D}}^{\prime}(\R)$}\,.
\end{equation}
One of the  main achievements of the paper \cite{DLR1} is the rewriting of  the equation \rec{ELE} in a more ``tractable" way in order to be able to investigate regularity and compactness property of 
weak  $1/2$-harmonic maps.
Precisely in \cite{DLR1} the following two results have been proved:

\begin{Proposition}
\label{pr-I.1}
A map $u$ in $\dot{H}^{1/2}(\R,{\cal{S}}^{m-1})$ is a weak $1/2$-harmonic map if and only if it satisfies the following   equation  
\begin{equation}\label{eq1intr}
(-\Delta)^{1/4} ( u\wedge (-\Delta)^{1/4}u)=T(u\wedge,u)\,,
\end{equation}
where, in general for arbitrary positive integers $n,m,\ell$, for every $Q\in \dot{H}^{1/2}(\R^n, {\cal{M}}_{\ell\times m}(\R))$ $\ell\ge 1$\footnote{ ${\cal{M}}_{\ell\times m}(\R)$ denotes, as usual, the space of $\ell\times m$ real matrices.} and $u\in\dot{H}^{1/2}(\R^n,\R^m)$ , $T$ is the operator defined by
\begin{equation}\label{opT}
T(Q,u):=(-\Delta)^{1/4}[Q\,(-\Delta)^{1/4}u]-Q(-\Delta)^{1/2} u+ (-\Delta)^{1/4} Q\,(-\Delta)^{1/4} u\,.\end{equation}
\hfill $\Box$
\end{Proposition}

The  equation (\ref{eq1intr}) has  been  completed by the following ''structure equation'' which is  a consequence of the fact that $u\in{{\cal{S}}^{m-1}}$ almost everywhere:

\begin{Proposition}
\label{pr-I.2}
All maps in $\dot{H}^{1/2}(\R,{\cal{S}}^{m-1})$ satisfy the following identity
\begin{equation}\label{eq2intr}
(-\Delta)^{1/4} (u\cdot  (-\Delta)^{1/4}u)=S(u\cdot,u)-{\cal{R}}((-\Delta)^{1/4} u \cdot {\cal{R}}(-\Delta)^{1/4} u)\,.
\end{equation}
where,  in general for arbitrary positive integers $n,m,\ell$, for every  $Q\in \dot{H}^{1/2}(\R^n, {\cal{M}}_{\ell\times m}(\R))$, $\ell\ge 1$ and $u\in\dot{H}^{1/2}(\R^n,\R^m)$, $S$ is the operator given by
\begin{equation}
\label{opS}
S(Q,u):=(-\Delta)^{1/4}[Q(-\Delta)^{1/4} u]-{\cal{R}} (Q\nabla u)+{\cal{R}}[(-\Delta)^{1/4} Q\,{\cal{R}}(-\Delta)^{1/4} u ]
\end{equation}
and   ${\cal{R}}$ is  the Fourier multiplier of symbol $m(\xi)=i\frac{\xi}{|{\xi}|}\,$.
\hfill $\Box$
\end{Proposition}

We call the operators $T$, $S$   {\em three-terms commutators} and  in \cite{DLR1} the following estimates have been established: for every $u\in \dot{H}^{1/2}({\R},{\R}^m)$ and 
$ Q\in H^{1/2}(\R,{\cal{M}}_{\ell\times m}(\R))$ we have
\begin{equation}
\label{zz3}
\|T(Q,u)\|_{H^{-1/2}(\R)}\le C\ \|Q\|_{\dot{H}^{1/2}(\R)}\ \|u\|_{\dot{H}^{1/2}(\R)}\,,
\end{equation}
\begin{equation}
\label{zz4}
\|S(Q,u)\|_{H^{-1/2}(\R)}\le C\ \|Q\|_{\dot{H}^{1/2}(\R)}\ \|u\|_{\dot{H}^{1/2}(\R)}\,,
\end{equation}
and
\begin{equation}
\label{zz5}
\|{\cal{R}}((-\Delta)^{1/4} u\cdot {\cal{R}}(-\Delta)^{1/4} u))\|_{\dot{H}^{-1/2}(\R)}\le C\  \|u\|_{\dot{H}^{1/2}(\R)}^2\,.
\end{equation}

\par 
 
 What has been discovered is a sort of 
   ''gain of regularity'' in the r.h.s of the equations \rec{eq1intr} and \rec{eq2intr} in the sense that, under the  assumptions $u\in \dot{H}^{1/2}({\R},{\R}^m)$ and 
$Q\in \dot{H}^{1/2}(\R,{\cal{M}}_{\ell\times m}(\R))$ each term individually in $T$ and $S$ - like for instance $(-\Delta)^{1/4}[Q(-\Delta)^{1/4}u]$ or $Q(-\Delta)^{1/2} u$ ... - is  not in $H^{-1/2}$
but the special linear combination of them constituting $T$ and $S$ is in $H^{-1/2}$.
The same phenomenon appears in 
  \underbar{in dimension 2} in the context of harmonic maps, for the Jacobians   $J(a,b):={\frac{\partial a}{\partial x} \frac{\partial b}{\partial y} -\frac{\partial a}{\partial y} \frac{\partial b}{\partial x}}$ (with $a,b\in \dot{H}^1(\R^2)$)  which satisfy  as a direct consequence of Wente's theorem (see \cite{CRW,W}
)\begin{equation}
\label{zz6}
\|J(a,b)\|_{\dot{H}^{-1}(\R^2)}\le C\ \|a\|_{\dot{H}^1(\R^2)}\ \|b\|_{\dot{H}^1(\R^2)}
\end{equation}
whereas, individually, the terms $\frac{\partial a}{\partial x} \frac{\partial b}{\partial y}$ and $\frac{\partial a}{\partial y} \frac{\partial b}{\partial x}$ are not in $H^{-1}(\R^2)$.

\medskip

The estimates \rec{zz3} and \rec{zz4} imply in particular that if $u\in\dot{H}^{1/2}(\R,{\cal{S}}^{m-1})$ is a $1/2$-harmonic map then
\begin{equation}\label{bt}
\|(-\Delta)^{1/4}u\|_{L^{2}(\R)}\le C \|(-\Delta)^{1/4}u\|^2_{L^{2}(\R)}\,.
\end{equation}
where the constant $C$ is independent on $u$.\par
From the inequality \rec{bt} 
it follows that if $ C \|(-\Delta)^{1/4}u\|_{L^{2}(\R)}<1$ then the solution is constant. This the so-called {\em bootstrap test} and it is the key observation to prove Morrey-type estimates and to 
deduce H\"older regularity of $1/2$-harmonic maps, see \cite{DLR1}\,.
 \par
 We mention here that since the paper \cite{DLR1} several extensions have been considered. The regularity of solutions to nonlocal linear Schr\"odinger systems with applications to $1/2$-harmonic maps with values into general manifolds have been studied by Rivi\`ere and the author in \cite{DLR2}.   $n/2$-harmonic maps in odd  dimension $n$ has been
 considered in \cite{Sch1}  and \cite{DL} respectively in the case of values into the $m-1$ dimensional sphere and into general manifolds and the case of $\alpha$-harmonic maps in $\dot W^{\alpha,p}(\R^n,{\cal{S}}^{m-1})$, with
 $\alpha p=n$\,, has been recently studied by Schikorra and  the author in \cite{DLS}\,.
 Finally Schikorra \cite{Sch2} has also studied the partial regularity of weak solutions to nonlocal linear systems with an antisymmetric potential in the supercritical case (namely where $\alpha p<n)$ under a crucial monotonicity assumption  on the solutions which allows to reduce to the critical case.\par
 In this paper we address to the issue  of  understanding the behaviour of sequences
 $u_k$ of weak $1/2$-harmonic maps.
 We observe that as in the case of harmonic maps the bootstrap test \rec{bt} implies that if the energy is {\em small} then the system behaves locally like a linear system of the form
 $(-\Delta)^{1/2} u=0$ (namely the r.h.s is ``dominated " by the l.h.s of the equation)\,.
 As a consequence we obtain that any sequence  $u_k$ of weak $1/2$-harmonic maps with uniformly bounded energy
 weakly converges to a weak $1/2$-harmonic map  $u_{\infty}$ and strongly converges to
  $u_{\infty}$ away from a finite (possibly empty)  set  $\{a_1 ,\ldots,a_{\ell} \}\subset\R$\,.\par
  Namely we have 
  ( up to a subsequence) 
$$(|(-\Delta)^{1/4}u_k|^2\rightharpoonup|(-\Delta)^{1/4}u_{\infty}|^2)dx+\sum_{i=1}^{\ell}\lambda_i \delta_{a_i},~~\mbox{in Radon measure}\,,$$
as $k\to+\infty$,  with $\lambda_i\ge 0\,.$
  It remains the question to understand how the convergence at the concentration points 
  $a_i $ fails to be strong\,. A careful analysis shows that the loss of energy during the weak convergence is not only concentrated at the points  $a_i $ but it is also quantized : this amount of energy is given by the sum of energies of non-constant $1/2$-harmonic maps (the so-called {\em bubbles})\,.
  More precisely we get the following result

\begin{Theorem}\label{main} 
Let $u_k\in \dot H^{1/2}(\R,{\cal{S}}^{m-1})$ be a sequence of $1/2$-harmonic maps such that 
$\|u_k\|_{\dot H^{1/2}}\le C$\,. Then it holds:\par
\begin{enumerate}
\item There exist $u_{\infty}\in \dot H^{1/2}(\R,{\cal{S}}^{m-1})$ and a possibly empty set $\{a_1,\ldots, a_\ell\}$, $\ell\ge 1\,,$ such that up to subsequence
\begin{equation}\label{conv}
u_{n }\to u_\infty \quad \mbox{in $\dot W^{1/2,p}_{loc}(\R\setminus\{a_1,\ldots, a_{\ell}\}),~~p\ge 2$ as $k \to +\infty$}
\end{equation}
and 
\begin{equation}\label{provv1}
(-\Delta)^{1/2} u_{\infty} \wedge u_{\infty}=0,\quad \mbox{in ${\cal{D}}^{\prime}(\R )$}
\end{equation}
\item There is a family $ \tilde u_{\infty}^{i,j}\in \dot H^{1/2}(\R,{\cal{S}}^{m-1})$ of $1/2$-harmonic maps $ ( i\in\{1,\ldots,\ell\}, j\in\{1,\ldots ,N_i\}), $ such that up to subsequence
\begin{equation}\label{finalquantintr}
\|(-\Delta)^{1/4}(u_k-u_{\infty}-\sum_{i,j} \tilde u_{\infty}^{i,j})\|_{L^2_{loc}(\R)}\to 0,~~\mbox{as $k \to +\infty$}\,.
\end{equation}
\end{enumerate}
\end{Theorem}
Theorem \ref{main} says that for every $i$, $\lambda_i=\sum_{j=1}^{N_i}L(\tilde u_{\infty}^{i,j})\,.$ 
Therefore there is no  dissipation of energy in the region between $u_{\infty}$  and the bubbles and between the bubbles themselves (the so-called {\em neck-regions})\,.\par
 We would like now to mention a result obtained in 
   the    paper \cite{DLR3} on the characterization  of  $1/2$ harmonic maps $u\colon {\cal{S}}^1\to {\cal{S}}^2$
   which permits us to deduce that   in the case of $1/2$ harmonic maps with
   values in ${\cal{S}}^2$ one has    $\lambda_i=2\pi n_i$, with $n_i$ a positive integer    and also to provide a simple example  showing that the quantization may actually occur, namely the set $\{a_1,\ldots, a_\ell\}$ may be nonempty\,.

        \begin{Theorem}\label{mindisk}\cite{DLR3}
        i) $u\colon {\cal{S}}^1\to {\cal{S}}^1$ is a weak  $1/2$-harmonic map    if and if its harmonic extension $\tilde u\colon D^2 \to \R^2$ is holomorphic or anti-holomorphic\,.\par
ii)  $u\colon {\cal{S}}^1\to {\cal{S}}^2$ is a weak  $1/2$-harmonic map if and if the composition of weak  $1/2$-harmonic map $u\colon {\cal{S}}^1\to {\cal{S}}^1$   and an isometry $T\colon {\cal{S}}^2\to {\cal{S}}^2$\,.  ~\hfill$\Box$
 \end{Theorem}
 We remark that because of the invariance of the Lagrangian \rec{nonloclagr} with respect to the trace of conformal transformations we can  study  without restrictions the problem in ${\cal{S}}^1$ instead of $\R$\,.\par 
From Theorem  \ref{mindisk} it follows that  $1/2$-harmonic maps $ u\colon {\cal{S}}^1\to {\cal{S}}^1$ with  $deg(u)=1$ coincide with the trace of M\"obius transformations of the disk $D^2\subseteq\R^2\,.$ Moreover  every non-constant  weak  $1/2$-harmonic map  $u\colon {\cal{S}}^1\to {\cal{S}}^2$        satisfies
 $$
 \int_{{\cal{S}}^1}|(-\Delta)^{1/4}u|^2 dx =2\pi k<+\infty\,.$$
where $k $ is a positive integer which coincides with $|deg(u)|$\,.

    Let us consider now the following sequence of $1/2$-harmonic maps
 $$u_n\colon {\cal{S}}^1\to {\cal{S}}^1,~~u_n(z)=\frac{z-a_n}{1-\bar{a_n}z},$$ with $|a_n|=1$ and $a_n\to 1$ as $n\to+\infty\,.$ 
 In this case we have $u_n\to -1$ in $C^{\infty}_{loc}(\R\setminus\{1\})\,,$ thus the set of concentration points is  nonempty. Theorem \ref{main} yields  the existence of one Bubble $\tilde u_\infty$ such that
 $$
 \|(-\Delta)^{1/4}(u_k-\tilde u_{\infty})\|_{L^2_{loc}}\to 0,~~~\mbox{as $n\to +\infty\,.$}$$
 
 We explain now the method we have used to prove the main Theorem \ref{main}\,.\par
 In order to get the quantization of the energy we exploit a ``functional analysis" method introduced by Lin and Rivi\`ere in\cite{LR} in the context of harmonic maps in no-conformal dimensions.
Such a method consists in the use of the interpolation Lorentz spaces in the special case where the r.h.s of the equation can be written as linear combination of Jacobians\,.\par

This techniques has been recently applied in \cite{LR1,LR2}  and in \cite{BR} for the quantization analysis   respectively of  linear Schr\"odinger
systems with antisymmetric potential in $2$-dimension,  of bi-harmonic maps in $4$-dimensions.
and of  Willmore surfaces. We refer the reader to the papers \cite{Riv,LR} for an overview of the bubbling and quantization issues in the  literature\,.\par
We describe briefly the key steps to get the quantization analysis.\par
1. First of all we will make use of a general result proved in \cite{BR} which permits to
split the domain (in our case $\R$)  into
{\em the converging region} (which is the complement of small neighborhoods of the $a^i$), {\em bubbles domains} and {\em neck-regions}  (which are unions of degenerate annuli).\par
2. We prove that the $L^2$ norm of $(-\Delta)^{1/4}u_k$  in the neck regions is arbitrary small, (see Theorem \ref{neck1})\,.   
Thanks to the duality of the Lorentz spaces $L^{2,1}-L^{2,\infty}$\footnote{The $L^{2,\infty}(\R)$ is the space of measurable functions $f$ such that
$$
\sup_{\lambda>0}\lambda |\{x\in\R~: |f(x)|\ge\lambda\}|^{1/2}<+\infty\,.
$$
$L^{2,1}(\R)$ is the Lorentz space of measurable functions satisfying
$$\int_{0}^{+\infty}|\{x\in\R~: |f(x)|\ge\lambda\}|^{1/2} d\lambda<+\infty \,\,.$$
}
 this is reduced in estimating the $L^{2,\infty}, L^{2,1}$ norms of $u_k$ in these regions.
Precisely we first show that the $L^{2,\infty}$ norm of the $u_k$ is arbitrary small in degenerate annuli (see Lemma \ref{neck2}) and
as far as the  $ L^{2,1}$ norm is concerned, we use  
  the following improved estimate on the operators $T$ and $S$  which is proved in the Appendix\,
 
 \begin{Theorem}\label{comm1}
Let $u,Q\in H^{1/2}((\R^n))$ . Then $T(Q,u),S(Q,u)\in {\cal{H}}^{1}(\R^n)$\footnote{${\mathcal H}^1({\R}^n)$ denotes   the Hardy space which is the space of $L^1$
functions $f$ on ${\R}^n$satisfying 
\[
\int_{{\R}^n}\sup_{t\in {\R}}|\phi_t\ast f|(x)\ dx<+\infty\quad ,
\]
where $\phi_t(x):=t^{-n}\ \phi(t^{-1}x)$ and where $\phi$ is some function in the Schwartz space ${\mathcal S}({\R}^n)$ satisfying $\int_{{\R}^n}\phi(x)\ dx=1$.\par
 {For more properties on the Hardy space ${\mathcal H}^1$ we refer to \cite{Gra1} and \cite{Gra2}.} }
 and
\begin{equation}\label{commest2}
\|T(Q,u)\|_{{\cal{H}}^{1}(\R^n)}\le C\|Q\|_{\dot{H}^{1/2}(\R^n)}\|u\|_{\dot{H}^{1/2}(\R^n)}\,.\end{equation}
\begin{equation}\label{commest3}
\|S(Q,u)\|_{{\cal{H}}^{1}(\R^n)}\le C\|Q\|_{\dot{H}^{1/2}(\R^n)}\|u\|_{\dot{H}^{1/2}(\R^n)}\,.~~~\hfill \Box
\end{equation}
\end{Theorem}

In a forthcoming paper \cite{DL3}  we are going to  investigate  bubbles and quantization issues in
the case  of nonlocal Schr\"odinger linear systems with applications to $1/2$-harmonic maps with values into manifolds\,. The difficulty there is to succeed in getting  an uniform  $L^{2,1}$ estimate as well on degenerate annuli as in the local case (see \cite{LR1}).
\par It would  be also very interesting   to understand the geometric properties of the bubbles in the case of more general manifolds\,.\par
This paper is organized as follows. \par
 
In Section 1 we address to the compactness issue which is  the first part of Theorem \ref{main}. In Section 2  we prove $L^2$ estimates on
degenerate annual domains. In Section 3  we prove the second part of Theorem \ref{main}.
In the Appendix we prove Theorem \ref{comm1}\,.

\par
\medskip
 
 \section{Compactness}
\par\bigskip
In this Section we prove the first part of Theorem \ref{main}\,. The result is based on the following {\em $\varepsilon$-regularity} property whose proof can be found in \cite{DLR1,DLR2} and in \cite{DL2}\,.
 
\begin{Lemma}[$\varepsilon$-regularity ]\label{epsreg}
Let $u\in \dot H^{1/2}(\R,{\cal{S}}^{m-1})$ be a $1/2$-harmonic map. Then there exists $\varepsilon_0>0$ such that if for some $\gamma>0$ 
$$\sum_{j\ge 0} 2^{-\gamma j}\|(-\Delta)^{1/4} u\|_{L^2(B(x,2^{j}r)) }\le \varepsilon_0\,,$$
then there is $p>2$ such that for every $x\in\R$, $y\in B(x,r/2)$ we have
\begin{equation}
\left(r^{p/2-1}\int_{B(y,r/2)}|(-\Delta)^{1/4} u(x)|^p dx\right)^{1/p}\le C\,\|u\|_{ \dot H^{1/2}(\R)}\,.
\end{equation}
\end{Lemma}
By bootstrapping into the equations \rec{eq1intr} and \rec{eq2intr} and by localizing Theorems \ref{comm2} and \ref{comm3} (see \cite{DL2}) one can show the following:
\begin{Corollary}
Let $u\in \dot H^{1/2}(\R,{\cal{S}}^{m-1})$ be a $1/2$-harmonic map. Then $(-\Delta)^{1/4} u$ is in $ L^p_{loc}(\R)\cap L_{loc}^{\infty}{(\R)}$ for every $p\ge 2$ with
\begin{equation}
r^{1/2-1/p}\|(-\Delta)^{1/4} u\|_{L^p(B(x,r)) }\le C\|u\|_{\dot{H}^{1/2}(\R)}\,,
\end{equation}
\begin{equation}
r^{1/2}\|(-\Delta)^{1/4} u\|_{L^{\infty}(B(x,r)) }\le C\|u\|_{\dot{H}^{1/2}(\R)}\,,
\end{equation}
for all $x\in\R$ and $r>0\,.$
\end{Corollary}

We will use  also the localized version of 
 the following result whose proof can be found in \cite{AH} page 78: 
\begin{Lemma}\label{AH}
Let $0<\alpha<1$ and $f\in L^p({\R})$, $1<p<\infty$. Then there is a constant $C>0$ independent on $f$, such that 
$$
\|\Delta^{-\frac{\alpha\theta}{2}} f\|_{L^r(\R^n)}\le C\, \|\Delta^{-\frac{\alpha}{2}} f\|_{L^s(\R^n)}^{\theta}\,\|f\|_{L^p(\R^n)}^{1-\theta}\,,$$
for $0<\theta<1$, $1\le s\le\infty$, $\frac{1}{r}=\frac{\theta}{s}+\frac{1-\theta}{p}\,.$
\end{Lemma}

 Next we show that if $\|u_k\|_{\dot H^{1/2}(B(x,\rho))}\le \varepsilon_0$ where $\varepsilon_0 >0$  is the small constant appearing in the $\varepsilon$-regularity  Lemma \ref{epsreg}, then
 $u_k\in \dot  W_{loc}^{\frac{1}{q}+\frac{1}{2},2}(\R)$\,.
\begin{Proposition}\label{epr} 
 Let $\rho>0$ be such that  $\|u_k\|_{\dot H^{1/2}(B(x,\rho))}\le \varepsilon_0$ with $\varepsilon_0 >0$ given in Lemma  \ref{epsreg}. Then for all $q>2$ there exists $C>0$ (independent on $x\,,k$) such that
 \begin{equation}\label{ests}
 \|(-\Delta)^{ \frac{1}{2q}+\frac{1}{4}} u_k\|_{L^2(B(x,\rho/2))}\le C\|u_k\|_{\dot{H}^{1/2}(\R)}\,.
 \end{equation}
 \end{Proposition}
  
 {\bf Proof of Proposition \ref{epr}\,.}
We set $v_k=(-\Delta)^{1/4}u_k$.  From Lemma \ref{epsreg} it follows that 
  there exists $q>2$ (independent on $k$) such that  for every $y\in B(x,\rho/2)$
\begin{equation}\label{unifest0}
\|  v_k\|_{L^{q}(B(y,\rho/4))}\le  C\|u_k\|_{\dot{H}^{1/2}(\R)}\,
\end{equation}
and bootstrapping into  the equations \rec{eq1intr} and \rec{eq2intr} one gets 
\begin{equation}\label{unifest}
\|(-\Delta)^{1/4} v_k\|_{L^{\frac{2q}{q+2}}(B(y,\rho/4))}\le C\|u_k\|_{\dot{H}^{1/2}(\R)}\,.
\end{equation}
   Now we set $f_k:=(-\Delta)^{1/4} v_k$. By applying Lemma \ref{AH}  in  $B(y,\rho/4)$ with $f:=v_k$ and $p=q$, $r=2$, $s=\frac{2q}{q+2}$, $\alpha=\frac{1}{2}$ and $\theta=\frac{q-2}{q}$
     we obtain
\begin{equation}\label{unifest2}
\|(-\Delta)^{-\frac{q-2}{4q}} f_k\|_{L^{2}(B(y,\rho/4))}\le C\|f_k\|_{L^{\frac{2q}{q+2}}(B(y,\rho/4))}\|v_k\|_{L^{q}(B(y,\rho/4))}
\end{equation}
  In particular we get that $v_k\in \dot W^{ \frac{1}{q},2}(B(y,\rho/4))$ and hence $u_k\in  W^{\frac{1}{q}+\frac{1}{2},2}(B(y,\rho/4))$
with
\begin{equation}\label{unifest3}
\| u_k\|_{\dot  W^{\frac{1}{q}+\frac{1}{2},2}(B(y,\rho/4))}=\| (-\Delta)^{1/2q+1/4}u_k\|_{L^{2}(B(y,\rho/4))}\le C\,.
\end{equation}
Actually one can show that  the estimate \rec{unifest3} holds for every $q>2\,.$ This concludes
the proof.~\hfill$\Box$
\par
\medskip

We  show now   a singular point removability type result for $1/2$-harmonic maps.\par\medskip
  \begin{Proposition}\label{removsing}[Singular point removability]
  Let $u \in \dot H^{ \frac{1}{2}} (\R,{\cal{S}}^{m-1})$ be a $1/2$-harmonic map in ${\cal{D}}^{\prime} (\R\setminus\{a_1,\ldots,a_{\ell}\})$. Then
  $$
   u \wedge  (-\Delta)^{ \frac{1}{2}}u =0 \quad \mbox{in ${\cal{D}}^{\prime} (\R)$\,.}$$
   \end{Proposition}
   {\bf Proof of Proposition \ref{removsing}\,.}
   The fact that
   $$
u \wedge  (-\Delta)^{ \frac{1}{2}}u =0 \quad \mbox{in ${\cal{D}}^{\prime} (\R\setminus\{a_1,\ldots,a_{\ell}\})$}$$ 
implies that
$$
 (-\Delta)^{1/4}(u \wedge  (-\Delta)^{1/4}u )=T(u \wedge ,u)~~\mbox{in ${\cal{D}}^{\prime} (\R\setminus\{a_1,\ldots,a_{\ell}\})$}\,,$$
 where
 $T(u \wedge,u)\in \dot H^{- \frac{1}{2}}(\R)$ and 
 \begin{equation}
 \|T(u \wedge,u)\|_{\dot H^{- \frac{1}{2}}(\R)}\le C \|u\|^2_{\dot H^{ \frac{1}{2}}(\R)}\,.
\end{equation}
The distribution
 $\phi:= (-\Delta)^{1/4}(u \wedge  (-\Delta)^{1/4}u)- (T(u \wedge,u))$    is of  order $p=1$ and supported in $\{a_1,\ldots,a_{\ell}\}$.
 Therefore by Schwartz Theorem \cite{Bony} one has
 $$\phi =\sum_{|\alpha|\le 1} c_{\alpha}\partial^{\alpha} \delta_{a_i}\,.
 $$
 Since $\phi\in \dot H^{- \frac{1}{2}}(\R) $,      then  the above implies that $c_{\alpha}=0$\footnote{Suppose by contradiction that a distribution $\phi\in H^{- \frac{1}{2}}(\R) $ satisfies
 $\phi=\sum_{|\alpha|\le 1}c_{\alpha}\partial^{\alpha}\delta_{a_0}$. We can write $\phi=(-\Delta)^{1/4}f$ for some $f\in L^2(\R)$. Then
 ${\cal{F}}[\phi](\xi)=|\xi|^{1/2}{\cal{F}}[f](\xi)=\sum_{|\alpha|\le 1}c_{\alpha}(i)^{|\alpha|}\xi^{\alpha} $ and this is not possible since ${\cal{F}}[f]\in L^2(\R)\,.$ Therefore $c_\alpha=0\,.$}
 is  and thus
 $$
 (-\Delta)^{1/4}(u\wedge (-\Delta)^{1/4}u )=T(u\wedge, u) ~~\mbox{in ${\cal{D}}^{\prime} (\R )$}\,.$$ 
 We conclude the proof of Proposition  \ref{removsing}\,.~\hfill$\Box$
 
The proof of the {\bf first part of Theorem \ref{main}} concerning the compactness of uniformly bounded $1/2$-harmonic maps is contained in the following Lemma\,.\par
\begin{Lemma}\label{sing} 
Let $u_k\in \dot H^{1/2}(\R,{\cal{S}}^{m-1})$ be a sequence of $1/2$-harmonic maps such that 
$\|u_k\|_{\dot H^{1/2}}\le C$\,. Then there exist a sequence $u_{k^{\prime}}$ of $u_k$, a function $u_{\infty}\in \dot H^{1/2}(\R,{\cal{S}}^{m-1})$ and $\{a_1,\ldots ,a_{\ell}\} $,
$\ell\ge 1$, such that
\begin{equation}\label{conv}
u_{k^{\prime}}\to u_\infty ~\mbox{as}~ k^{\prime}\to +\infty~  \mbox{in $\dot H^{\frac{1}{2}}_{loc}((\R\setminus\{a_1,\ldots, a_{\ell}\}))$ ~~for $p\ge 2$,  }
\end{equation}
and 
\begin{equation}\label{provv1}
(-\Delta)^{1/2} u_{\infty} \wedge u_{\infty}=0,\quad \mbox{in ${\cal{D}}^{\prime}(\R)\,.$}
\end{equation}
\end{Lemma}
{\bf Proof of Lemma \ref{sing}\,.}\par 
{\bf 1.}  First of all there exists a subsequence $u_{k^{\prime}}$ of $u_k$, a function $u_{\infty}\in \dot H^{1/2}(\R,{\cal{S}}^{m-1})$ such that $u_{k^{\prime}}\rightharpoonup u_\infty$ as $k^{\prime}\to +\infty\,.$\par
{\bf 2.} If $\|u_k\|_{\dot H^{1/2}(B(x,\rho))}\le \varepsilon_0$ for all $k\ge 1$  then from Lemma \ref{epsreg} and the  Rellich-Kondrachov Theorem (if  $\Omega\subset \R$ is a bounded subset then the embedding $W^{\frac{1}{q}+\frac{1}{2},2}(\Omega) \hookrightarrow  W^{ 1/2,t}(\Omega)$ is compact 
for all $t<\frac{2q}{q-2}$)  it follows that 
$$u_{k^{\prime}}\to u_{\infty} ~\mbox{as}~ k^{\prime}\to +\infty~\mbox{in}~ \dot H^{1/2}(B(x,\rho/4),{\cal{S}}^{m-1}).$$
for all $x\in\R$\,.  In particular we have
$$\mbox{$(-\Delta)^{ \frac{1}{2}} u_{k^{\prime}}\to (-\Delta)^{ \frac{1}{2}}u_{\infty}$ as $k^{\prime}\to +\infty$ in $\dot H^{-1/2}(B(x,\rho/4),{\cal{S}}^{m-1})$}\,.$$ Hence 
$$
(-\Delta)^{ \frac{1}{2}} u_{k^{\prime}}\wedge  u_{k^{\prime}}\to (-\Delta)^{ \frac{1}{2}} u_{\infty}\wedge u_\infty~\mbox{as}~ k^{\prime}\to +\infty~  \mbox{in   ${\cal{D}}^{\prime}(B(x,\rho/4))$\,,}
$$
and
$$
(-\Delta)^{ \frac{1}{2}} u_{\infty}  \wedge u_\infty=0 \quad \mbox{in   ${\cal{D}}^{\prime}(B(x,\rho/4))$}\,.
$$
{\bf 3.}  {\bf Claim 1:} There are only finitely many points $\{a_1,\ldots,a_\ell\}$ such that
\begin{equation}\label{uinfty}
(-\Delta)^{ \frac{1}{2}} u_{\infty} \wedge u_{\infty}=0,\quad \mbox{in ${\cal{D}}^{\prime}(\R\setminus\{a_1,\ldots, a_{\ell}\})$}
\end{equation}
{\bf Proof of the claim 1.} We associate to every $x$ the number $\rho_x^n>0$ such that
$\|u_k\|_{\dot H^{1/2}(B(x,\rho_x^n)}=\varepsilon_0$ where $\varepsilon_0$ is as in Lemma \ref{epsreg}\,.
\par
For every $M>0$ and $n\ge 1$ we set
$$I_k^M:=\{x~:~\rho_x^n<\frac{1}{M} \}$$
and
$$
{\cal{F}}_k^M:=\{B(x,\frac{1}{M}),~x\in I_k^M\}\,.$$
 By Vitali-Besicovitch Covering Theorem (see for instance \cite{EG}), we can find an at most countable family of points $(x_j^{k,M})_{j\in  J_k^M}$, $x_j^{k,M}\in  I_k^M$ and
  and $I_k^M\subseteq\cup_{j\in J_k^M}B(x_j^{k,M},\frac{1}{M})\,.$ Moreover every  $x\in  I_k^M$ is contained in at most $K$
   balls, $K$ being a number depending only on the dimension of the space\,. \par
   Now we observe that
   \begin{eqnarray*}
   C&\ge &\|u_k\|^2_{\dot H^{1/2}(\R)}\\
   &\ge& \sum_{j\in J_k^M}\iint_{B(x_j^{k,M},\frac{1}{M})\times B(x_j^{k,M},\frac{1}{M})}\frac{|u_k(x)-u_k(y)|^2}{|x-y|^2} dx dy\\
   &\ge & \sum_{j\in J_k^M}\varepsilon_0^2=|J_k^M|\varepsilon_0^2\,.
   \end{eqnarray*}
   Thus $|J_k^M|<+\infty$ for every $k$ and $M$ and this implies that for $k$ and $M$ large enough $|J_k^M|=C$, with $C$ independent on $k$ and $M\,.$ In particular there exists
   $k_0>0$ such that
   $$
   I_k^M\subseteq \cup_{j=1}^{k_0} B(x_j^{n,M},\frac{1}{M})\,.$$
   By definition we have $I_k^{M+1}\subseteq I_k^{M}$ for all $n$ and $M$. By using a {\em diagonal procedure} we can subtract 
   a subsequence $k^{\prime}\to +\infty$ such that $x_{j}^{k^{\prime},M}\to x_{j}^{\infty,M} $ for all $M>0$ and  $j$ and
   $$
   I_{\infty}^M\subseteq  \cup_{j=1}^{k_0} B(x_j^{\infty,M},\frac{1}{M})\,.$$
   Now we let $M\to +\infty$ and get
   $$
   I_{\infty}^0\subseteq  J_{\infty,0}:= \{x_j^{\infty,0}\}_{j=1,\ldots n_0}\,.$$
   {\bf Claim 2.} If $x\notin J_{\infty,0}$ then there exits $\tilde r>0$ such that
   $$
   u_\infty\wedge  (-\Delta)^{ \frac{1}{2}}u_{\infty}=0 \quad \mbox{in ${\cal{D}}^{\prime} (B(x,\tilde r))$\,.}$$
   {\bf Proof of the claim 2.}
   We assume that $x_{j}^{\infty,0}\ne \infty$ for all $j=1,\ldots,n_0\,.$
   Let $\gamma=dist(x,J_{\infty,0})$ and $K>0$ be such that $2K^{-1}<\gamma\,.$
   Let $\tilde M>0$ be such that for all $M\ge \tilde M$ and for all   $j=1,\ldots,n_0$ we have
\begin{equation}\label{xi0}
   |x_j^{\infty,0}-x_j^{\infty,M}|< \frac{1}{4K}
   \,.\end{equation}
   Let $\bar k>0$ be such that for all $k^{\prime}\ge \bar k$ and for all   $j=1,\ldots,n_0$ we have
   
    \begin{equation}\label{xiM}
   |x_j^{\infty, M}-x_j^{k^{\prime}, M}|< \frac{1}{4K}
     \,.\end{equation}
  
  By combining \rec{xi0} and \rec{xiM} we get
  \begin{eqnarray*}
     |x -x_j^{k^{\prime}, M}|&\ge & |x-x_j^{\infty,0}|
     \\
     && -  |x_j^{\infty,0}-x_j^{\infty,  M}|-  |x_j^{\infty,\tilde M}-x_j^{k^{\prime},  M}|\\
     &\ge & \frac{2}{K}-\frac{1}{4K}-\frac{1}{4K}=\frac{3}{2K}>\frac{1}{M}\,.
     \end{eqnarray*}
     Therefore $x\notin \cup_{j=1}^{n_0} B(x_j^{k^{\prime},  M},\frac{1}{  M})$,   and $x\notin I_{k^{\prime}}^{  M}$ for all $k\ge\bar k$\,.
     In particular $\rho_{k^{\prime},x}\ge {  M}^{-1}$ and (up to subsequence) $\rho_{k^{\prime},x}\to \rho_{\infty,x}>0$.
     Now let $0<\tilde r<\rho_{\infty,x}$. Then 
     $$
     B(x,\tilde r)\subseteq B(x,\rho_{k^{\prime},x}),\quad \mbox{for $k$ large}\,.$$
  Since we have $\|u_{k^{\prime}}\|_{\dot H^{1/2} B(x,\rho_{k^{\prime},x}))}=\varepsilon_0$, then by applying Step 2 we get
 $$
   u_\infty\wedge  (-\Delta)^{1/2}u_{\infty}=0 \quad \mbox{in ${\cal{D}}^{\prime} (B(x,\tilde r))$\,.}$$
 This concludes the proof of the claim 2   by setting $a_i:=x_i^{\infty}$ and $\ell=n_0$ \,.\par
 Now we apply Proposition \ref{removsing} and we get that
 $$
   u_\infty\wedge  (-\Delta)^{1/2}u_{\infty}=0 \quad \mbox{in ${\cal{D}}^{\prime} (\R)$\,.}$$
 We can conclude the proof of the Lemma \ref{sing} and the first part of Theorem \ref{main}\,.
  ~\hfill$\Box$\par\medskip
   
 \section{$L^{2,\infty}$ and $L^2$  estimates in degenerate annuli}
   In this Section we will prove some energy estimates of $1/2$-harmonic maps in degenerate annuli. Such estimates are crucial in the next Section in order to get the quantization analysis in the neck regions\,.\par
   The main result of this Section is
   \begin{Theorem}\label{neck1}
     There exists  $ \tilde\delta >0$ such that for any  $1/2$-harmonic maps $u\in \dot H^{ \frac{1}{2}}(\R,{\cal{S}}^{m-1})$, for any $\delta<\tilde\delta$  and $\lambda,\Lambda>0$ with $\lambda<(2\Lambda)^{-1} $  satisfying 
 \begin{equation}\label{quant0}
 \sup_{\rho\in[\lambda,(2\Lambda)^{-1} ]}\left( \int_{B(0,2\rho)\setminus B(0,\rho)}| (-\Delta)^{1/4}u|^2 dx\right)^{1/2}\le \delta \,,\end{equation}
    we have
 \begin{eqnarray}\label{quant2}
&&\int_{B(0,\Lambda^{-1})\setminus B(0,\lambda )}| (-\Delta)^{1/4}u|^2dx 
\\
&&~~~~~~\le C   \sup_{\rho\in[\lambda,(2\Lambda)^{-1} ]}\left( \int_{B(0,2\rho)\setminus B(0,\rho)}| (-\Delta)^{1/4}u|^2 dx\right)^{1/2}\,.\nonumber\end{eqnarray}
  \end{Theorem}
  The proof of Theorem \ref{neck1}  consists in three steps: \par
  1) first we show that we can control in degenerate annuli 
  the $L^q$ norm of $(-\Delta)^{1/4}u$ for some $q>2$ by
  \begin{equation}\label{estprovv}
  \sup_{\rho\in[\lambda,(2\Lambda)^{-1} ]}\left( \int_{B(0,2\rho)\setminus B(0,\rho)}| (-\Delta)^{1/4}u|^2 dx\right)^{1/2}\,.\end{equation}\par
   2) then we estimate the $L^{2,\infty}$ norm of
   $ (-\Delta)^{1/4}u$ in degenerate annuli in terms of \rec{estprovv},\par
    3)  finally we use the global $L^{2,1}$ estimates obtained in the appendix (see Theorem \ref{comm1}) and the duality $L^{2,1}-L^{2,\infty} in order to conclude\,.$

 \begin{Lemma}  [$L^q$-estimates]\label{qesti}
 There exists  $ \tilde\delta >0$ such that for any  $1/2$-harmonic maps $u\in \dot H^{ \frac{1}{2}}(\R)$, for any $\delta<\tilde\delta$, $\lambda,\Lambda>0$ with $2\lambda<(4\Lambda)^{-1} $ such that
 \begin{equation}\label{quant1}
 \sup_{\rho\in[\lambda,(2\Lambda)^{-1} ]}\left( \int_{B(0,2\rho)\setminus B(0,\rho)}| (-\Delta)^{1/4}u|^2 dx\right)^{1/2}\le \delta \,.\end{equation}
 then   there exists $q>2$ (independent on $\lambda,\Lambda,u$)  such that
 \begin{eqnarray*}\label{quant2}
 &&\sup_{\rho\in[2\lambda ,(4\Lambda)^{-1}]}\left(\rho^{q/2-1}\int_{B(0,2\rho)\setminus B(0,\rho)}| (-\Delta)^{1/4}u|^qdx\right)^{1/q}\\
 &&~~~~~~~~\le C  \sup_{\rho\in[\lambda,(2\Lambda)^{-1} ]}\left( \int_{B(0,2\rho)\setminus B(0,\rho)}| (-\Delta)^{1/4}u|^2 dx\right)^{1/2}\,.\nonumber\end{eqnarray*}

    \end{Lemma}
 {\bf Proof of Lemma \ref{qesti}\,.} We choose $\delta=\frac{\varepsilon_0}{2}$ where $\varepsilon_0>0$ is the constant appearing in the $\varepsilon$-regularity Lemma \ref{epsreg}\,. \par
 {\bf Step 1.} {There exist $p>2$ (independent on $\lambda,\Lambda,u$)   such that}
  \begin{equation}\label{quant3}
 \sup_{\rho\in[2\lambda ,(4\Lambda)^{-1}]}\left(\rho^{p/2-1}\int_{B(0,2\rho)\setminus B(0,\rho)}| (-\Delta)^{1/4}u|^pdx\right)^{1/p}\le C \|u\|_{\dot H^{ \frac{1}{2}}(\R)}\,.\end{equation}
  
 {\bf Proof of Step  1.}
 
 Let
 $r>0$ be such that
 $$
 \left( \int_{B(0,2r)\setminus B(0,r)}| (-\Delta)^{1/4}u|^2 dx\right)^{1/2}<\delta\,.$$
 {\bf Claim:} There exists $p>2$ (independent on $\delta$ and $r$) such that
 \begin{equation}\label{quant2}
 \left(r^{p/2-1}\int_{B(0,\frac{3}{2}r)\setminus B(0,\frac{5}{4}r)}| (-\Delta)^{1/4}u|^p dx\right)^{1/p}\le  C \|u\|_{\dot H^{ \frac{1}{2}}(\R)}\,.
 \end{equation}
Let $y\in B(0,\frac{3}{2}r)\setminus B(0,\frac{5}{4}r)$, (we clearly have   $dist(y,\partial (B(0,2r)\setminus B(0, r)))\ge 1/4$). \par
 Let $j_0\ge 3$ such that $2^{-j_0/2}\left(  \int_{\R}| (-\Delta)^{1/4}u|^2 dx\right)^{1/2}\le \delta\,,$ and
 $B(y,2^{-j_0}r)\subset (B(0,2r)\setminus B(0,r))$ for all $y\in B(0,\frac{3}{2}r)\setminus B(0,\frac{5}{4}r)$\,. \par\bigskip
{\bf Estimate of $\sum_{h\ge 0}2^{-h/2}\|(-\Delta)^{1/4}u\|_{L^2(B(0,2^h(2^{-2j_0}r))}$\,.}
 \begin{eqnarray*}
 &&\sum_{h\ge 0}2^{-h/2}\|(-\Delta)^{1/4}u\|_{L^2(B(0,2^h(2^{-2j_0}r))}\\
&=& \sum_{h=0}^{ j_0}2^{-h/2}\|(-\Delta)^{1/4}u\|_{L^2(B(0,2^h(2^{-2j_0}r))}+\sum_{h= j_0+1}^{\infty}2^{-h/2}\|(-\Delta)^{1/4}u\|_{L^2(B(0,2^h(2^{-2j_0}r))}\\
&\leq&\delta ( \sum_{h=0}^{ \infty}2^{-h/2} )
+2^{-(j_0+1)/2} \left( \int_{\R}| (-\Delta)^{1/4}u|^2 dx\right)^{1/2}\\
&\le & \delta +\delta =2\delta =\varepsilon_0\,.
\end{eqnarray*}
Now we apply Lemma \ref{epsreg} : there exists $p>2$ and $C_{j_0}>0$  such that
\begin{equation}\label{quant3}
 \left(r^{p/2-1}\int_{B(y,2^{-(j_0+1)}r) }| (-\Delta)^{1/4}u|^p dx\right)^{1/p}\le  C_{j_0} \|u\|_{\dot H^{ \frac{1}{2}}(\R)}\,.
 \end{equation}
 By covering the annulus $B(0,\frac{3}{2}r)\setminus B(0,\frac{5}{4}r)$ by a finite number of balls $B(y,2^{-(j_0+1)}r)$ we finally get
 $$
 \left(r^{p/2-1}\int_{B(0,\frac{3}{2}r)\setminus B(0,\frac{5}{4}r)}| (-\Delta)^{1/4}u|^p dx\right)^{1/p}\le  \tilde C_{j_0} \|u\|^2_{\dot H^{ \frac{1}{2}}(\R)}\,,
$$
and {\bf the proof of Claim 1} is concluded.\par\medskip

Hence  
 \begin{equation}\label{quant4}
 \sup_{\rho\in[2\lambda ,(4\Lambda)^{-1}]}\left(\rho^{p/2-1}\int_{B(0,2\rho)\setminus B(0,\rho)}| (-\Delta)^{1/4}u|^pdx\right)^{1/p}\le C \|u\|_{\dot H^{ \frac{1}{2}}(\R)}\,.\end{equation}
 We thus conclude the {\bf proof of Step  1} \,.\par
 
 {\bf Step 2.} There exists $q>2$ (independent on $\lambda,\Lambda,u$ and dependent on $p$)  such that
 \begin{eqnarray}\label{quant2bis}
&& \sup_{\rho\in[2\lambda ,(4\Lambda)^{-1}]}\left(\rho^{q/2-1}\int_{B(0,2\rho)\setminus B(0,\rho)}| (-\Delta)^{1/4}u|^qdx\right)^{1/q}\\
&&~~~\le C  \sup_{\rho\in[\lambda,(2\Lambda)^{-1} ]}\left( \int_{B(0,2\rho)\setminus B(0,\rho)}| (-\Delta)^{1/4}u|^2 dx\right)^{1/2}\,.\nonumber\end{eqnarray}

  {\bf Proof of Step 2\,.}
   Let us take $q^{-1}=\theta p^{-1} +(1-\theta)2^{-1}$.
 Then by H\"older Inequality  and by using \rec{quant3} we get
 \begin{eqnarray*}
&&\left(\rho^{q/2-1}\int_{B(0,2\rho)\setminus B(0,\rho)}| (-\Delta)^{1/4}u|^pdx\right)^{1/q}\\
&\le &
 \left(\rho^{p/2-1}\int_{B(0,2\rho)\setminus B(0,\rho)}| (-\Delta)^{1/4}u|^pdx\right)^{1/p}
 \left(\int_{B(0,2\rho)\setminus B(0,\rho)}| (-\Delta)^{1/4}u|^2dx\right)^{1/2}\\
 &\le& C \|u\|_{\dot H^{ \frac{1}{2}}(\R)}  \sup_{\rho\in[\lambda,(2\Lambda)^{-1} ]}\left( \int_{B(0,2\rho)\setminus B(0,\rho)}| (-\Delta)^{1/4}u|^2 dx\right)^{1/2}\\
 &\le& C    \sup_{\rho\in[\lambda,(2\Lambda)^{-1} ]}\left( \int_{B(0,2\rho)\setminus B(0,\rho)}| (-\Delta)^{1/4}u|^2 dx\right)^{1/2}\,.
 \end{eqnarray*}
 This  concludes the {\bf proof of Step 2 and of  Lemma \ref{qesti}}\,.~\hfill$\Box$
 \begin{Lemma}[$L^{2,\infty}$ estimates]\label{neck2}
   There exists  $ \tilde\delta >0$ such that for any  $1/2$-harmonic maps $u\in \dot H^{ \frac{1}{2}}(\R)$, for any $\delta<\tilde\delta$  and $\lambda,\Lambda>0$  with $\lambda<(2\Lambda)^{-1} $ satisfying
 $$ \sup_{\rho\in[\lambda,{(2\Lambda)}^{-1}]}\left( \int_{B(0,2\rho)\setminus B(0,\rho)}| 
 (-\Delta)^{1/4} u|^2dx\right)^{1/2}\le \delta $$
then
\begin{equation}\label{linftyest}
 \| (-\Delta)^{1/4} u\|_{L^{2,\infty}(B(0,(2\Lambda)^{-1})\setminus B(0,\lambda )}\le  C \sup_{\rho\in[\lambda,{(2\Lambda)}^{-1}]}\left( \int_{B(0,2\rho)\setminus B(0,\rho)}| 
 (-\Delta)^{1/4} u|^2dx\right)^{1/2}  \,.
\end{equation}
 where $C$ is independent on $\rho,u,\lambda,\Lambda\,.$
 \end{Lemma}
 {\bf Proof of Lemma \ref{neck2}.}
 We set $  f=(-\Delta)^{1/4} u$ in $B(0,(4\Lambda)^{-1})\setminus B(0,2\lambda )$ and $ f=0$ otherwise.\par
 
  Let $ \delta<\tilde\delta/4$ where $\tilde\delta$ is the constant appearing in Theorem \ref{neck1}\,.  From Lemma  \ref{qesti}    it  follows that 
 for all $\lambda,\Lambda>0$ with $2\lambda<(4\Lambda)^{-1}$ if 
  $$ \sup_{\rho\in[\lambda,{(2\Lambda)}^{-1}]}\left( \int_{B(0,2\rho)\setminus B(0,\rho)}| 
 (-\Delta)^{1/4} u|^2dx\right)^{1/2}\le \delta\, $$
then   there exists
$q>2$,  such that
 \begin{eqnarray*}
 &&  \sup_{\rho\in[2\lambda ,(4\Lambda)^{-1}]}\left(\rho^{q/2-1}\int_{B(0,2\rho)\setminus B(0,\rho)}| 
 (-\Delta)^{1/4} u|^qdx\right)^{1/q} \\
 &&~~~~~\le C \sup_{\rho\in[\lambda,{(2\Lambda)}^{-1}]}\left( \int_{B(0,2\rho)\setminus B(0,\rho)}| 
 (-\Delta)^{1/4} u|^2dx\right)^{1/2} \,. 
    \end{eqnarray*}
We set $$\gamma=C \sup_{\rho\in[2\lambda,{(4\Lambda)}^{-1}]}\left( \int_{B(0,2\rho)\setminus B(0,\rho)}| 
 (-\Delta)^{1/4} u|^2dx\right)^{1/2}$$
 We observe that   for all $\rho\in [2\lambda ,(4\Lambda)^{-1}]$  one has:
 \begin{eqnarray*}
\gamma^q&\ge & \rho^{q/2-1} \int_{B(0,2\rho)\setminus B(0,\rho)}| f|^qdx\\
&\ge&\rho^{q/2-1} \alpha^q|\{x\in B(0,2\rho)\setminus B(0,\rho):~|f|>\alpha\}|\,.
\end{eqnarray*}
Let $k\in\Z$, then the following estimate holds

\begin{eqnarray*}
&&\alpha^2\sum_{j\ge k}|\{x\in B(0,2^{j+1}\alpha^{-2})\setminus B(0,2^j\alpha^{-2}):~|f|>\alpha\}|\\[5mm]
& &\le
\alpha^2\sum_{j\ge k}\frac{(2^{j+1}\alpha^{-2})^{1-q/2}}{\alpha^q}\gamma^q
 = \gamma^q\sum_{j\ge k} 2^{(j+1)(1-q/2)}\\[5mm]
 &&\le  \gamma^q 2^{1-q/2} 2^{k(1-q/2)}(1-2^{1-q/2})^{-1}\\[5mm]
 && =\gamma^q2^{k(1-q/2)}(2^{q/2-1}-1)^{-1}\,.
 \end{eqnarray*}
 Therefore
 \begin{eqnarray*}
\alpha^2 |\{x\in \R:~|  f|>\alpha\}|&\le& \gamma^q 2^{1-q/2} 2^{k(1-q/2)}+\alpha^2|B(0,2^k\alpha{-2})|\\[5mm]
 &\le & \gamma^q 2^{k(1-q/2)}(2^{q/2-1}-1)^{-1}+\alpha^2 2 2^{k} \alpha^{-2}\,.
 \end{eqnarray*}
 Now we choose $k$ in such a way that $2^{k}=\gamma^2/2$. It follows that
 $$
 \alpha^2 |\{x\in \R:~|  f|>\alpha\}|\le \frac{ \gamma^2}{2}(2^{q/2-1}-1)^{-1}  + {\gamma^{2}} =
 \frac{2^{q/2}-1}{2^{q/2}-2}\gamma^2\,.
  $$
 Hence
 \begin{eqnarray}\label{linftyest2}
 & &\| (-\Delta)^{1/4} u\|_{L^{2,\infty}(B(0,(4\Lambda)^{-1})\setminus B(0,2\lambda  )}\nonumber\\[5mm]
 &&~~=\sup_{\alpha>0} (\alpha^2 |\{x\in B(0,(4\Lambda)^{-1})\setminus B(0,2\lambda):~| (-\Delta)^{1/4} u(x)  |>\alpha\})^{1/2}\\[5mm]
& & ~~\le  \left(\frac{2^{q/2}-1}{2^{q/2}-2}\right)^{1/2}\gamma\,.\nonumber
 \end{eqnarray}
By combining \rec{linftyest2} and the fact that the $L^{\infty}$ norms of $(-\Delta)^{1/4} u $ in the annuli $B(0,\Lambda^{-1})\setminus B(0,(4\Lambda^{-1}))$ and $B(0,2\lambda )\setminus B(0,\lambda)$ are controlled by the
respective $L^2$ norms we get the estimate \rec{linftyest} and we conclude the proof of the Lemma \ref{neck2}\,.~\hfill$\Box$
 \par\medskip
 Now we can prove Theorem \ref{neck1}.\par
 {\bf Proof of Theorem \ref{neck1}\,.}
 Form Theorem \ref{comm1} it follows that any $1/2$-harmonic map $u\in \dot H^{1/2}(\R,{\cal{S}}^{m-1})$ satisfies $\|(-\Delta)^{1/4}u\|_{L^{2,1}(\R)}\le C$ where
 $C$ depends on $\|u\|_{\dot H^{1/2}(\R,{\cal{S}}^{m-1})}\,.$
 \par
 Now it is enough to use the duality $L^{2,1}-L^{2,\infty}$ and Lemma \ref{neck2} to get
 \begin{eqnarray*}
&& \int_{B(0,(\Lambda)^{-1})\setminus B(0,\lambda)}| (-\Delta)^{1/4}u|^2 dx\le 
 \|u\|_{L^{2,1}(\R)} \| (-\Delta)^{1/4} u\|_{L^{2,\infty}(B(0,(2\Lambda)^{-1})\setminus B(0,\lambda  )}\nonumber\\
&&~~~\le C \sup_{\rho\in[\lambda,(2\Lambda)^{-1} ]}\left( \int_{B(0,2\rho)\setminus B(0,\rho)}| (-\Delta)^{1/4}u|^2 dx\right)^{1/2}\,.\nonumber
 \end{eqnarray*}
 We can conclude  the proof of Theorem \ref{main}\,.~\hfill $\Box$
 
 
 \section{Bubbles and neck-regions }
   \par

   In the proof of the first part of Theorem \ref{main} (see Lemma \ref{sing}) we have shown   ( up to a subsequence) that  
$$|(-\Delta)^{1/4}u_k|^2\rightharpoonup|(-\Delta)^{1/4}u_{\infty}|^2dx+\sum_{i=1}^{\ell}\lambda_i \delta_{a_i},~~\mbox{in Radon measure}\,.$$
  The aim of this Section is to  
  show that   for every $i\in\{1\ldots\ell \}$ there exist bubbles $( \tilde {u}_{\infty}^{i,j}),$
$j\in\{1,\ldots,N_i\}$ such that $\lambda_i=\sum_{j=1}^{N_i}\int_{\R}|(-\Delta)^{1/4} \tilde {u}_{\infty}^{i,j}|^2 dx\,.$\par
We first give   the following  definitions.
 \begin{Definition}[Bubble]
 A {\bf Bubble} is a {\bf non-constant}    $1/2$-harmonic map $u\in \dot H^{1/2}(\R,{\cal{S}}^{m-1})$\,.   \end{Definition}
 \begin{Definition}[Neck region]
 A neck region for a function $f\in L^2(\R)$ is the union of finite   degenerate annuli of the type $A_k(x)=B(x,R_k)\setminus B(x,r_k)$ with $r_k\to 0$ and $\frac{R_k}{r_k}\to +\infty$ as $k \to +\infty$ satisfying the following property:  for all $\delta>0$ there exists  $\Lambda>0$ such that 
 $$
\left( \sup_{\rho\in[\Lambda r_k,(2\Lambda)^{-1}R_k]}\int_{B(x,2\rho)\setminus B(x,\rho)}|f|^2 dx\right)^{1/2}\le \delta\,.
 $$
 \end{Definition}

 \par
 \medskip

 {\bf Proof of the second part of Theorem \ref{main} \,.}\par
We have to show that
 there is a family $ \tilde u_{\infty}^{i,j}\in \dot H^{1/2}(\R,{\cal{S}}^{m-1})$, of non-constant $1/2$-harmonic maps $ ( i\in\{1,\ldots,\ell\}, j\in\{1,\ldots ,N_i\}), $ such that up to subsequence
\begin{equation}\label{finalquant}
\|(-\Delta)^{1/4}(u_k-u_{\infty}-\sum_{i,j} \tilde u_{\infty}^{i,j})\|_{L^2_{loc}(\R)}\to 0,~~\mbox{as $n\to\infty$}\,.
\end{equation}
We first observe  that
the {\em bootstrap test} \rec{bt} implies that each bubble has a bounded from below energy $c_0>0$\,. Therefore for every $i$,  $N_i<+\infty.$

 For simplicity we assume that $\ell=1$ and that there are at most two bubbles.\par

Now let us take  $\delta<\tilde\delta$ such that $C\delta<\varepsilon_0$ (the constants $C$ and $\tilde\delta$ are the one appearing in the statement of Theorem \ref{neck1})\,.\par
We also set $\gamma=\min(\frac{\delta}{2},\frac{\varepsilon_0}{2})\,.$\par
{\bf Step 1.} 
For every $n\ge 1$ we set
$$
 \rho^{1}_{k}=\inf\{\rho>0:~\exists x\in B(a_1,1): ~\int_{B(x,\rho)}|(-\Delta)^{1/4}u_k|^2dx =\gamma\}$$
\underline{There are two cases:}\par
{\bf Case 1.:} $\liminf_{k \to +\infty}  \rho^{1}_{k} >0$\par
In this case there is not concentration of the energy, namely  $\lambda_1=0$.\par
{\bf Case 2.:} $\lim_{k \to +\infty}  \rho^{1}_{k}=0$. For every $k\ge 1$, let $x^1_k\in B(a_1,1)$ be the point such that
$\int_{B( x_{1,k}, \rho^{1}_{k})}|(-\Delta)^{1/4}u_k|^2dx =\gamma\,.$ We have (up to subsequence)
$x^1_k\to a_1$ as $k\to+\infty\,$ (outside any neighborhood of $a_1$ there is no concentration)\,.\par
Now we choose a subsequence of $u_k$ (that we still denote by $u_k$) and a fixed radius $\alpha>0$ such that
$$
\limsup_{n\to\infty}\left[\sup_{0<r<\alpha}\{\int_{B(a_1,\alpha)\setminus B(a_1,r)}|(-\Delta)^{1/4} u_k(y)|^2dy  =\gamma\}\right]=0\,.
$$
Now we borrow the idea in \cite{BR}   to  split the annulus $B(x_{1,k},\alpha)\setminus B(x_{1,k},\rho_{k}^1)$ in domains of unbounded conformal class where the energy is small and domains of bounded conformal class where the energy is bounded from below.\par

Precisely by applying Lemma 3.2 in \cite{BR}, we can find    a sequence of family of radii
 $$
 R_k^0=\alpha>R_k^1>\ldots>R^{N_1}_k=\rho^1_{k}\,$$
 with $\{1,\ldots,N_1\}=I_0\cup I_1\,.$
 For every $i_\ell\in I_0$ one has
 \begin{equation}
 \lim_{k \to +\infty}\log\left( \frac{{R_k}^{i_\ell}}{{R_k}^{i_\ell+1}}\right)<+\infty
 ~~\mbox{and}~~\int_{B(x_{1,k},{R_k}^{i_\ell+1})\setminus B(x_{1,k}, {R_k}^{i_\ell})}|(-\Delta)^{1/4}\tilde u_k(y)|^2dy 
 \ge \gamma\,,
 \end{equation}
and for every $i_\ell\in I_1, $ one has
 \begin{equation}
 \lim_{k \to +\infty}\log \left(\frac{{R_k}^{i_\ell}}{{R_k}^{i_\ell+1}}\right)=+\infty
 ~~\mbox{and}~~\forall\rho\in({R_k}^{i_\ell},{R_k}^{i_\ell+1}/2),~
 \int_{B(x_{1,k},2\rho )\setminus B(x_{1,k}, \rho)}|(-\Delta)^{1/4}\tilde u_k(y)|^2dy \le 2\gamma\,.
 \end{equation}
 We consider the smallest annulus $A_k^{i_\ell}:=B(x_{1,k},{R_k}^{i_\ell})\setminus B(x_{1,k}, {R_k}^{i_\ell+1})$ of the first type $i_\ell\in I_0$. For such an $i_\ell$ we define 
 $$
r^{i_\ell}_{k}=\inf\{r<{R_k}^{i_\ell+1}:~\exists x\in A_k^{i_\ell}: ~\int_{B(x,r)}|(-\Delta)^{1/4}u_k|^2dx =\gamma\}\,.$$
We consider the following two cases.\par
{\bf Case 1.}  There exists a subsequence of $ r^{i_\ell}_{k}$ such that
$$
\lim_{k\to+\infty}\frac{r^{i_\ell}_{k}}{{R_k}^{i_\ell}}>0\,.
$$
In this case there is not concentration of the energy  in $A_k^{i_\ell}$ and we pass to the next 
$A_k^{i^{\prime}_\ell}$ (if there is any).\par
{\bf Case 2.}   We have
$$
\lim_{k\to+\infty}\frac{r^{i_\ell}_{k}}{{R_k}^{i_\ell}}=0\,.
$$
In this case we have once again concentration. Let $x_{2,k}\in A_k^{i_\ell}$ such that
$$
\int_{B(x_{2,k},r^{i_\ell}_{k})}|(-\Delta)^{1/4}u_k|^2dx =\gamma\,.$$

and we  set $\rho_{k}^{2}=r^{i_\ell}_{k}$\,.  
\par

{\bf We separate two sub-cases:}\par
 \par
{\bf \underline{Case of  two ``separated" bubbles}}   $\liminf_{k \to +\infty} \frac{\rho_{k}^{1}}{\rho_{k}^{2}}>0\,.$
In this case the following two conditions hold
$$
\left\{\begin{array}{cc}
\lim_{n\to\infty}\frac{|x_{1,k}-x_{2,k}|}{\rho_k^1 }=+\infty\\
[5mm]
\lim_{n\to\infty}\frac{|x_{1,k}-x_{2,k}|}{\rho^2_k}=+\infty\,.\end{array}
\right.
$$
In this case  the bubbles $\tilde u_{2,\infty}$ and $\tilde u_{1,\infty}$ are ``independent" \,.\par
\medskip
  \par
  Let us consider the 
  two ``separated" balls $B(x_{1,k},\rho_{k}^{1})$
 and $B(x_{2,k},\rho_{k}^{2})$, with
 $$
 \lim_{k\to\infty}\frac{|x_{1,k}-x_{2,k}|}{\rho_{k}^{1}+\rho_{k}^{2}}=+\infty\,.$$
 For every $\alpha$ we set 
 $$
 {\cal{N}}^1_k(\alpha)=B(a_1,\alpha)\setminus \left(B(x_{1,k},\alpha^{-1}\rho_{k}^{1 })\cup B(x_{2,k},\alpha^{-1}\rho_{k}^{ 2})\right)\,.$$
 The above construction gives the existence of $\alpha$ small enough independent of $k$  such that 
 $$
\left\{\begin{array}{c}
 \mbox{for} j=1,2 ~~\mbox{and for all $\rho$ such that $B(x_{j,k},2\rho)\setminus 
 B(x_{j,k},2\rho)\subseteq  {\cal{N}}^1_k(\alpha)$}\\[5mm]
 \int_{B(x_{j,k},2\rho)\setminus 
 B(x_{j,k},\rho)}|\Delta)^{1/4}u_k|^2dx\le  2\gamma\,
  \\[5mm]
 ~~\mbox{and}\\[5mm]
  \int_{B(x_{j,k},\rho_{k}^{j})}|\Delta)^{1/4}u_k|^2dx=\gamma\,.
  \end{array}
  \right.$$
  {\bf Claim:} 
 the region $ {\cal{N}}^1_k(\alpha)$ is a  {neck-region}\,. 
 \par
  {\bf Proof of the Claim: } it  is a consequence of the  following general property.\par
  \begin{Lemma}\label{neck3}
   Let $A_k=B(x_k,R_k)\setminus B(x_k,r_k)$ an annulus satisfying
 $r_k\to 0$, $\frac{R_k}{r_k}\to +\infty$ and $x_k\to x_{\infty}$ as $k \to +\infty\,,$ and
 \begin{equation}\label{gammacond}
  \sup_{r_k\le\rho\le\frac{R_k}{2}}\int_{B(x_{ k},2\rho)\setminus 
 B(x_{k},\rho)}|(-\Delta)^{1/4}u_k|^2dx\le  2\gamma\,.\end{equation}
 Then for all $\eta>0$ there exists $\Lambda>0$ such that
\begin{equation}\label{etacond3}
\sup_{r \in [ \Lambda r_k,\Lambda^{-1}R_k]}\int_{B( x_{ n},2r )\setminus B( x_{ n}, r)}|(-\Delta)^{1/4}u_k|^2dx\le  \eta\,.
\end{equation}
\end{Lemma}
{\bf Proof of Lemma \ref{neck3}\,.}
  Suppose by contradiction that there exists $\eta >0$ and two sequence $\Lambda_k\to +\infty$  as $k\to+\infty$ and $\Lambda_k r_k\le \tilde r_{ k}\le (\Lambda_k)^{-1} R_k $     such that
\begin{equation}\label{etacond2}
\int_{B( x_{k},2\tilde r_{k} )\setminus B( x_{ n}, \tilde r_{k})}|(-\Delta)^{1/4}u_k|^2dx>\eta\,
\end{equation}
We define $\tilde u_{k}(y)=u(\tilde r_{k} y+ x_{1,k})$\,.
From  condition \rec{gammacond}
 and Theorem \ref{neck1} it follows that
   \begin{equation} 
   \int_{B(0,\frac{R_k}{\tilde r_k})\setminus 
 B(0,\frac{r_k}{\tilde r_k})}|(-\Delta)^{1/4}\tilde u_k|^2dx\le  2\gamma\,.\end{equation}
 Lemma \ref{epsreg} and Lemma \ref{sing} imply that 
$ \tilde u_{k}\to \tilde u_{\infty}$ in $\dot W^{1/2,p}_{loc}(\R\setminus \{0\})$ for all $p\ge 1\,$, where
 $\tilde u_{\infty} $ is a nontrivial $1/2$-harmonic maps ($\int_{B( 0,2  )\setminus B( 0, 1)}|(-\Delta)^{1/4}\tilde u_{\infty} |^2dx>\eta$)\,. On the other  hand the condition    \rec{gammacond} gives    $$
 \int_{\R}|(-\Delta)^{1/4}\tilde u_\infty|^2dx\le C\delta< \varepsilon_0\,.$$
 The bootstrap test yields that $\tilde u_\infty$ is trivial which is a contradiction.

We conclude the proof of Lemma \ref{neck3} and of the claim.~\hfill $\Box$\par
\medskip
By applying Theorem \ref{main}  we get that for all $\eta>0$ small enough
$$
\int_{{\cal{N}}^1_k(\alpha)}|(-\Delta)^{1/4} u_k^2dx\le \eta\,.
$$

 \medskip
{\bf\underline{ Case of bubble over bubble}\,.}   $\liminf_{k \to +\infty} \frac{\rho_{k}^{1}}{\rho_{k}^{2}}=0\,.$
We define $\tilde u_{2,k}(y)=u(\rho_{k}^2 y+ x_{2,k})$ We have 
$\tilde u_{2,k}\to \tilde u_{2,\infty}$ in $W^{1/2,p}_{loc}(\R\setminus \{a_{1}\})$ for all $p\ge 1\,$ and  
$\int_{B(0,2)\setminus B(0,1)}|(-\Delta)^{1/4}\tilde u_{2,\infty}|^2dx\ge \delta$. Therefore $\tilde u_{2,\infty}$
is a new bubble   
(case {\em Bubble over Bubble}: the bubble $\tilde u_{2,\infty}$ ``contains " $\tilde u_{1,\infty}$ and for $k$ large enough $x_{1,k}\in B(x_{1,2}, \rho_{k}^{2 })$)\,.
For every $\alpha$ we set 
  $$
 {\cal{N}}^{1,2}_k(\alpha)= \left(B(x_{1,k},\alpha\rho_{k}^{2 })\setminus B(x_{1,k},\alpha^{-1}\rho_{k}^{ 1})\right)$$
and 
$$ {\cal{N}}_k(\alpha)= {\cal{N}}^{1}_k(\alpha)\cup  {\cal{N}}^{1,2}_k(\alpha)\,.$$
By arguing as above one can show that $ {\cal{N}}_k(\alpha)$ is a neck region.

Since we have assumed that there are at most two bubbles, the procedure stops here.
 Otherwise one has to continue the procedure until annuli of the type $I_0$ have been explored.

 Therefore  for every $\eta>0$ we get \par
  \underline{1. case of independent bubbles:}
 \begin{eqnarray}\label{fin1}
&&\lim_{k\to+\infty} \int_{\R}|(-\Delta)^{1/4}  u_k|^2dx=
\lim_{k\to+\infty} \int_{{\cal{N}}^1_k(\alpha)}|(-\Delta)^{1/4} u_k|^2dx\nonumber\\&&+
\sum_{j=1}^2\lim_{k\to+\infty} \ \int_{B(x_{j,k},\alpha^{-1}\rho_{k}^{j})}|\Delta)^{1/4}u_k|^2dx\nonumber\\
&&+\lim_{k\to+\infty} \int_{\R\setminus B(a_1,\alpha)}|(-\Delta)^{1/4} u_k|^2dx\\
&&\le  \eta+\sum_{j=1}^2 \int_{B(0,\alpha^{-1})}|\Delta)^{1/4}\tilde u^{j}_\infty|^2dx\nonumber\\&&+\int_{\R\setminus B(a_1,\alpha)}|(-\Delta)^{1/4}  u_\infty|^2dx\,.\nonumber
\end{eqnarray} 
  \underline{2.  case of  bubble over bubble}
\begin{eqnarray}\label{fin2}
&&\lim_{k\to+\infty} \int_{\R}|(-\Delta)^{1/4}  u_k|^2dx\le 
\lim_{k\to+\infty} \int_{{\cal{N}}_k(\alpha)}|(-\Delta)^{1/4} u_k|^2dx\\&&+
\lim_{k\to+\infty} [\int_{B(x_{2,k},\alpha^{-1}\rho_{k}^{2})\setminus B(x_{1,k},\alpha \rho_{k}^{2})}|(-\Delta)^{1/4}u_k|^2dx+\int_{B(x_{1,k},\alpha^{-1}\rho_{k}^{1})}|(-\Delta)^{1/4}u_k|^2dx]
\nonumber\\
&&+\lim_{k\to+\infty} \int_{\R\setminus B(a_1,\alpha)}|(-\Delta)^{1/4}  u_k|^2dx\nonumber\\
&&\le  \eta+  \int_{B(0,\alpha^{-1})\setminus B(0,\alpha)}|(-\Delta)^{1/4}\tilde u^{2}_\infty|^2dx\nonumber\\&&+
 \int_{B(0,\alpha^{-1})} |(-\Delta)^{1/4}\tilde u^{1}_\infty|^2dx+\int_{\R\setminus B(a_1,\alpha)}|(-\Delta)^{1/4}  u_\infty|^2dx\,.\nonumber
\end{eqnarray}

By taking in \rec{fin1} and \rec{fin2} the $\lim$ for $\alpha,\eta\to 0$ we get   the desired quantization estimate
\rec{finalquant}\,. This concludes the proof of the second part of Theorem \ref{main}~\hfill$\Box\,.$\par

\appendix
\section{Commutator estimates: Proof of Theorem \ref{comm1}}
In this Section we prove   Theorem \ref{comm1}. To this end 
 we shall make use of the Littlewood-Paley  dyadic  decomposition of unity that we recall here. Such a decomposition can be obtained as follows \,.
 Let $\phi(\xi)$ be a radial Schwartz function supported in $\{\xi\in{\R}^n:~|\xi|\le 2\}$, which is
 equal to $1$ in $\{\xi\in{\R}^n: ~|\xi|\le 1\}$\,.
 Let $\psi(\xi)$ be the function given by
 $$
 \psi(\xi):=\phi(\xi)-\phi(2\xi)\,.
 $$
 $\psi$ is then a ''bump function'' supported in the annulus $\{\xi\in{\R}^n:~1/2\le |\xi|\le 2\}\,.$
 
 \medskip
 
 Let $\psi_0=\phi$, $\psi_j(\xi)=\psi(2^{-j}\xi)$ for $j\ne 0 \,.$ The functions $\psi_j$, for $j\in\Z$, are supported in  $\{\xi\in{\R}^n:~2^{j-1}\le |\xi|\le 2^{j+1}\}\,$
and they realize  a dyadic decomposition of  the unity :
  $$
  \sum_{j\in\Z}\psi_j(x)=1\,.
  $$
 We further denote   $$\phi_j(\xi):=\sum_{k=-\infty}^j\psi_k(\xi)\,.
 $$ 
 The function $\phi_j$ is supported on  $\{\xi, ~|\xi|\le 2^{j+1}\}$.\par

  For every $j\in \Z$ and $f\in{\cal{S}}^{\prime}(\R)$  we define the Littlewood-Paley projection operators $P_j$ and $P_{\le j}$ by 
   \begin{eqnarray*}
 \widehat{ P_jf}=\psi_j \hat{f}~~~\widehat{ P_{\le j}f}=\phi_j \hat{f}\,.
 \end{eqnarray*}
 Informally $P_j$ is a frequency projection to the annulus $\{2^{j-1}\le |\xi|\le 2^{j}\}$, while
 $P_{\le j}$ is a frequency projection to the ball $\{|\xi|\le 2^j\}\,.$ We will set
 $f_j=P_j f$ and $f^j=P_{\le j} f$\,.
 
 We observe that  $f^j=\sum_{k=-\infty}^{j} f_k$  and $f=\sum_{k=-\infty}^{+\infty}f_k$ (where the convergence is in ${\cal{S}}^\prime(\R)$)\,.\par
 Given $f,g\in {\cal{S}}^\prime(\R)$  we can  split  the product in the following way
 \begin{equation}\label{decompbis}
 fg=\Pi_1(f,g)+\Pi_2(f,g)+\Pi_3 (f,g),\end{equation}
 where
 \begin{eqnarray*}
 \Pi_1(f,g)&=& \sum_{-\infty}^{+\infty} f_j\sum_{k\le j-4} g_k= \sum_{-\infty}^{+\infty} f_j g^{j-4}\,;\\
 \Pi_2(f,g)&=& \sum_{-\infty}^{+\infty} f_j\sum_{k\ge j+4} g_k =\sum_{-\infty}^{+\infty} g_j f^{j-4}\,;\\
 \Pi_3(f,g)&=& \sum_{-\infty}^{+\infty}  f_j\sum_{|k-j|< 4} g_k\,.\
 \end{eqnarray*}
 We observe that   for every $j$ we have 
 $$\mbox{supp${\cal{F}}[f^{j-4}g_j]\subset \{2^{j-2}\le |\xi|\le 2^{j+2}\}$};$$
 $$\mbox{supp${\cal{F}}[\sum_{k=j-3}^{j+3}f_jg_k]\subset \{|\xi|\le 2^{j+5}\}$}\,.$$
 The three pieces of the decomposition \rec{decompbis} are examples of paraproducts. Informally the
 first paraproduct $\Pi_1$ is an operator which allows high frequences of $f$ $(\thicksim 2^j)$ multiplied by  low frequences of $g$ $(\ll 2^j)$ to produce high frequences in the output. The second paraproduct
 $\Pi_2$ multiplies low fequences of $f$ with high frequences of $g$ to produce high fequences in the output. The third paraproduct $\Pi_3$ multiply high frequences of $f$ with high frequences of $g$ to produce comparable or lower frequences in the output. For a presentation of these paraproducts we refer to the reader for instance  to the book \cite{Gra2}\,.
 The following two Lemmae will be often  used   in the sequel.  
 
\begin{Lemma}\label{max}
For every $f\in {\cal{S}}^{\prime}$ we have
 $$ \sup_{j\in Z}|f^j|\le M(f)\,.$$
 \end{Lemma}
 
  \begin{Lemma}\label{lemmaprel1}
Let    $\psi$ be a Schwartz  radial function such that $supp (\psi)\subset B(0,4)$. Then for every $s\ge \left[\frac{n}{2}\right]+1$  we have 
$$
\|(-\Delta)^{s} {\cal{F}}^{-1}\psi\|_{L^1}\le C_{\psi,n}(1+s^{n+1}) 4^{2s}\,,
$$
where $C_{\psi,n}$ is a positive constant depending on the $C^2$ norm of $\psi$ and the dimension\,.
\end{Lemma}
 
\par
\medskip
\begin{Lemma}\label{lemmaprel2}
Let $f\in B^0_{\infty,\infty}(\R^n)$\footnote{ The homogeneous Besov space $\dot{B}_{\infty,\infty}^0(\R^n)$ is the space of tempered distribution $u$ for which $ \|u\|_{\dot{B}_{\infty,\infty}^0(\R^n)};=\sup_{j\in \Z} \|{\cal{F}}^{-1}[\psi_j {\cal{F}}[u]]\|_{L^{\infty}(\R^n)}$ is finite, (see for the precise definition of the Besov spaces \cite{RS})}. Then for all  $s\ge \left[\frac{n}{2}\right]+1$  and for all $j\in Z$ we have 
$$
2^{-2s j}\|(-\Delta)^{s} f_j\|_{L^\infty}\le C_{\psi,n} (1+s^{n+1})    4^{2s} \|f\|_{ B^0_{\infty,\infty}(\R^n)}\,.
$$
\end{Lemma}
\par
 For the proof of Lemma \ref{max} we refer to \cite{DLR1} and of Lemmae \ref{lemmaprel1} and
 \ref{lemmaprel2} we refer to \cite{DL}\,.

Given $u,Q$  we introduce the following pseudodifferential operators
\begin{equation}\label{opT}
T(Q,u):=(- \Delta)^{1/4}(Q(- \Delta)^{1/4}u)-Q(- \Delta)^{1/2} u+(- \Delta)^{1/4} u(- \Delta)^{1/4} Q\,\end{equation}
and
\begin{equation}
\label{opS}
S(Q,u):=(-\Delta)^{1/4}[Q(-\Delta)^{1/4} u]-{\cal{R}} (Q\nabla u)+{\cal{R}}((-\Delta)^{1/4} Q{\cal{R}}(-\Delta)^{1/4} u )
\end{equation}
and   ${\cal{R}}$ is  the Fourier multiplier of symbol $m(\xi)=i\frac{\xi}{|{\xi}|}\,$.
We prove in this Section some estimates on the operators \rec{opT} and \ref{opS}\,.\par\medskip

  {\bf Proof of Theorem \ref{comm1}.}\par
We make the proof for $n=1$\,. The case $n>1$ is analogous (for the details we refer to \cite{DL2})\par
{$\bullet$}  Estimate of $\|\Pi_1((- \Delta)^{1/4}(Q(- \Delta)^{1/4} u)\|_{{\cal{H}}^1}$\,.\par
\begin{eqnarray}\label{commestpi1}
&&\|\Pi_1((- \Delta)^{1/4}(Q(- \Delta)^{1/4} u)\|_{{\cal{H}}^1}=\int_{\R^n}(\sum_{j=0}^{\infty} 2^j Q^2_j((- \Delta)^{1/4} u^{j-4})^2)^{1/2} dx\\
&&\le \int_{\R^n} \sup_{j}|(- \Delta)^{1/4} u^{j-4}|(\sum_j 2^j Q_j^2)^{1/2} dx\\
&&\le (\int_{\R^n} (M((- \Delta)^{1/4} u))^2 dx)^{1/2}(\int_R\sum_j 2^j Q_j^2 dx)^{1/2}\nonumber\\
&&\le C \|Q\|_{\dot{H}^{1/2}(\R)}\|u\|_{\dot{H}^{1/2}(\R)}\,.\nonumber
\end{eqnarray}
{$\bullet$}  Estimate of $\|\Pi_1( (- \Delta)^{1/4} Q (- \Delta)^{1/4} u )\|_{{\cal{H}}^1}$.  
\begin{eqnarray}\label{commestpi12}
&&\|\Pi_1((- \Delta)^{1/4}Q(- \Delta)^{1/4} u)\|_{{\cal{H}}^1}=\int_{\R^n}(\sum_{j=0}^{\infty} ((- \Delta)^{1/4} Q_j)^2((- \Delta)^{1/4} u^{j-4})^2)^{1/2} dx\\
&&\le \int_{\R^n} \sup_{j}|(- \Delta)^{1/4} u^{j-4}|(\sum_j ((- \Delta)^{1/4} Q_j)^2)^{1/2} dx\\
&&\le (\int_{\R^n} (M((- \Delta)^{1/4} u))^2 dx)^{1/2}(\int_R\sum_j ((- \Delta)^{1/4} Q_j)^2dx)^{1/2}\nonumber\\
&&\le C \|Q\|_{\dot{H}^{1/2}(\R )}\|u\|_{\dot{H}^{1/2}(\R )}\,.\nonumber
\end{eqnarray}

{$\bullet$}  Estimate of $\|\Pi_2( (- \Delta)^{1/4} Q (- \Delta)^{1/4} u )\|_{{\cal{H}}^1}$.  It is as in \rec{commestpi12}\,.
\par
 $\bullet$ Estimate of $\Pi_3((- \Delta)^{1/4}(Q(- \Delta)^{1/4} u))$\,.\par
We show that it is in $B^0_{1,1}$. 
We observe that if $h\in B^0_{\infty,\infty}$ then $(- \Delta)^{1/4} h\in  B^{-1/2}_{\infty\infty}$ and thus

\begin{eqnarray*}
&&\sup _{\|h\|_{B^0_{\infty,\infty}}\le 1}\int_{\R }\sum_j\sum_{|k-j|\le 3}(- \Delta)^{1/4}(Q_j(- \Delta)^{1/4} u_k)h\\
&&=\sup _{\|h\|_{B^0_{\infty,\infty}}\le 1}\int_{\R}\sum_j\sum_{|k-j|\le 3}(- \Delta)^{1/4}(Q_j(- \Delta)^{1/4} u_k)
[(- \Delta)^{1/4} h^{j-6}+\sum_{t= j-5}^{j+6}(- \Delta)^{1/4}h_t]dx
\end{eqnarray*}
We have
\begin{eqnarray*}
&&\sup _{\|h\|_{B^0_{\infty,\infty}}\le 1}\int_{\R^n}\sum_j\sum_{|k-j|\le 3}(- \Delta)^{1/4}(Q_j(- \Delta)^{1/4} u_k)
h^{j-6}dx\\
&& =\sup _{\|h\|_{B^0_{\infty,\infty}}\le 1}\int_{\R^n}\sum_j\sum_{|k-j|\le 3} (Q_j(- \Delta)^{1/4} u_k)
(- \Delta)^{1/4} h^{j-6}dx\\
&&\le C \sup _{\|h\|_{B^0_{\infty,\infty}}\le 1}\|h\|_{B^0_{\infty,\infty}}\int_{\R^n}\sum_j\sum_{|k-j|\le 3} 2^{j/2} Q_j(- \Delta)^{1/4} u_k dx\\
&&\le C(\int_{\R^n}\sum_j  2^{j} Q^2_j dx)^{1/2}(\int_{\R^n}\sum_j ((- \Delta)^{1/4} u_j)^2dx)^{1/2}\\
&&
\le C \|Q\|_{\dot{H}^{1/2}(\R)}\|u\|_{\dot{H}^{1/2}(\R)}\,.
\end{eqnarray*}
By analogous computations we get 
\begin{eqnarray*}
&&\sup _{\|h\|_{B^0_{\infty,\infty}}\le 1}\int_{\R^n}\sum_j\sum_{|k-j|\le 3}(- \Delta)^{1/4}(Q_j(- \Delta)^{1/4} u_k)
[\sum_{t= j-5}^{j+6}(- \Delta)^{1/4}h_t]dx\le C \|Q\|_{\dot{H}^{1/2}(\R)}\|u\|_{\dot{H}^{1/2}(\R)}\,.\end{eqnarray*}


$\bullet$ Estimate of $ \Pi_3((- \Delta)^{1/4}Q(- \Delta)^{1/4}  u)-Q(- \Delta)^{1/2}u)\,.$
\begin{eqnarray}\label{pi3crochet1}
&&\|\Pi_3((- \Delta)^{1/4}Q(- \Delta)^{1/4} u-Q(- \Delta)^{1/2}u)\|_{B^0_{1,1}}\\
&&=\sup _{\|h\|_{B^0_{\infty,\infty}}\le 1} 
\int_{\R^n}\sum_j\sum _{|k-j|\le 3}[(- \Delta)^{1/4}(Q_j(- \Delta)^{1/4} u_k)-Q_j(- \Delta)^{1/2}u_k)[h^{j-6}+\sum_{t=j-5}^{j+6}h_t]dx\nonumber
\end{eqnarray}
We only estimate the terms with $h^{j-6}$, being the estimates with $ h_t$ similar\,.
\begin{eqnarray}\label{pi3crochet2}
&& \sup _{\|h\|_{B^0_{\infty,\infty}}\le 1} \int_{\R^n}\sum_j\sum _{|k-j|\le 3}[(- \Delta)^{1/4}(Q_j(- \Delta)^{1/4} u_k)-Q_j(- \Delta)^{1/2}u_k)[h^{j-6}]dx\\&&=
\sup _{\|h\|_{B^0_{\infty,\infty}}\le 1} \int_{\R^n}\sum_j\sum _{|k-j|\le 3} {\cal{F}}[h^{j-6}]
{\cal{F}}[(-\Delta)^{1/4}Q_j (- \Delta)^{1/4} u_k-Q_j(- \Delta)^{1/2}u_k]dx\nonumber\\ && 
 =
 \sup _{\|h\|_{B^0_{\infty,\infty}}\le 1} \int_{\R^n}\sum_j\sum _{|k-j|\le 3} {\cal{F}}[h^{j-6}]\nonumber\\
 &&\left[\int_{\R^n}({\cal{F}}[Q_j](y){\cal{F}}[(- \Delta)^{1/4}u_k](x-y)(|y|^{1/2}-|x-y|^{1/2})dy\right]dx\,.\nonumber
 \end{eqnarray}
 Now we observe that in \rec{pi3crochet2} we have $|x|\le 2^{j-3}$ and $2^{j-2}\le |y|\le 2^{j+2}$.
Thus $|\displaystyle\frac{x}{y}|\le \frac{1}{2}\,.$ 
\par
Hence
\begin{eqnarray}\label{estkern}
|y|^{1/2}-|\  x-y|^{1/2}&=&|y|^{1/2}[1-|1-\frac{x}{y}|^{1/2}]\\
&=&|y|^{1/2}\frac{x}{y}[1+|1-\frac{x}{y}|^{1/2}]^{-1}\nonumber\\
&=&|y|^{1/2}\sum_{k=0}^\infty\frac{c_k}{k!}(\frac{ x}{y})^{k+1}\,.\nonumber
\end{eqnarray}
 We may suppose that $\sum_{k=0}^\infty\frac{c_k}{k!}(\frac{ x}{y})^{k+1}$ is convergent if $|\displaystyle\frac{x}{y}|\le \frac{1}{2}\,,$ otherwise one may consider a
 different Littlewood-Paley decomposition by replacing the exponent $j-4$ with $j-s$, $s>0$ large enough. 
 We    introduce the following notation: for every $k\ge 0$ we set
$$S_k g={\cal{F}}^{-1}[\xi^{-(k+1)}|\xi|^{1/2} {\cal{F}} g].$$
We note that  if $h\in B^s_{\infty,\infty}$ then
$S_k h\in B^{s+1/2+k}_{\infty,\infty}$ and if $h\in H^{s}$ then $S_k h\in H^{s+1/2+k}\,.$\par
Moreover if $Q\in H^{1/2}$ then $\nabla^{k+1}(Q)\in H^{-k-1/2}\,.$\par
 We continue the estimate \rec{pi3crochet2}\,.\par
 
\begin{eqnarray*}\label{pi3crochet3}
\rec{pi3crochet2}&&= \sup _{\|h\|_{B^0_{\infty,\infty}}\le 1} \int_{\R^n}\sum_j\sum _{|k-j|\le 3} {\cal{F}}[h^{j-6}]\nonumber\\
 &&\left[\int_{\R^n}({\cal{F}}[Q_j](y){\cal{F}}[(- \Delta)^{1/4}u_k](x-y)(|y|^{1/2}-|x-y|^{1/2})dy\right]dx
\nonumber\\
&& =
\sup _{\|h\|_{B^0_{\infty,\infty}}\le 1} 
\sum_{\ell=0}^\infty\frac{c_{\ell}}{\ell!}\int_{\R^n}
\sum_j\sum _{|k-j|\le 3}
(-i)^{\ell+1}{\cal{F}}[\nabla^{\ell+1}h^{j-6}]{\cal{F}}[S_{\ell}Q_j(- \Delta)^{1/4} u_k)]( x) d x\nonumber\\
&&
\le \sup _{\|h\|_{B^0_{\infty,\infty}}\le 1} 
\sum_{\ell=0}^\infty\frac{c_{\ell}}{\ell!}\int_{\R^n}
\sum_j\sum _{|k-j|\le 3}[\nabla^{\ell+1}h^{j-6} [S_{\ell}Q_j(- \Delta)^{1/4} u_k)]( x) d x\nonumber\\
&&\mbox{by Lemma \ref{lemmaprel2}}\nonumber\\
&&
\le C
\sup _{\|h\|_{B^0_{\infty,\infty}}\le 1} 
\sum_{\ell=0}^\infty\frac{c_{\ell}}{\ell!}2^{-6\ell}4^{\ell+1}\|h\|_{B_{\infty,\infty}^0}\nonumber\\
&&
\int_{\R^n}
\sum_j\sum _{|k-j|\le 3}2^{(\ell+1)j} [S_{\ell}Q_j(- \Delta)^{1/4} u_k)]( x) d x\nonumber\\
&&
\le
C \sum_{\ell=0}^\infty\frac{c_{\ell}}{\ell!}2^{-6\ell}4^{\ell+1}(\int_{\R^n}\sum_j 2^{2(\ell+1)j}|S_\ell Q_j|^2 dx)^{1/2}(\int_{\R^n}\sum_j |(- \Delta)^{1/4}u_j|^2 dx)^{1/2}\nonumber\\
&&\mbox{by Plancherel Theorem} \nonumber\\
&&
\le
C \sum_{\ell=0}^\infty\frac{c_{\ell}}{\ell!}2^{-6\ell}4^{\ell+1}(\int_{\R^n}\sum_j 2^{2(\ell+1)j}|{\cal{F}}[S_\ell Q_j]|^2 dx)^{1/2}(\int_{\R^n}\sum_j |(- \Delta)^{1/4}u_j|^2 dx)^{1/2}\\
&&
\le
C \sum_{\ell=0}^\infty\frac{c_{\ell}}{\ell!}2^{-6\ell}4^{\ell+1}(\int_{\R^n}\sum_j 2^{2(\ell+1)j}2^{2(1-j)(\ell+1/2)}|{\cal{F}}[Q_j]|^2 dx)^{1/2}(\int_{\R^n}\sum_j |(- \Delta)^{1/4}u_j|^2 dx)^{1/2}\nonumber\\
&&
\le
C \sum_{\ell=0}^\infty\frac{c_{\ell}}{\ell!}2^{-6\ell}4^{\ell+1}2^{\ell}(\int_{\R^n}\sum_j 2^{j} Q_j^2 dx)^{1/2}(\int_{\R^n}\sum_j |(- \Delta)^{1/4}u_j|^2 dx)^{1/2}\nonumber\\
&&\le C
  \sum_{\ell=0}^\infty\frac{c_{\ell}}{\ell!}2^{-3\ell} \|Q\|_{\dot{H}^{1/2}(\R )}\|u\|_{\dot{H}^{1/2}(\R )}  \nonumber
 \end{eqnarray*}

\medskip
$\bullet$ Estimate of $\Pi_2((- \Delta)^{1/4}(Q(- \Delta)^{1/4} u)-Q(- \Delta)^{1/2}u))\,.$

\begin{eqnarray}\label{pi2crochet1}
&&\|\Pi_2((- \Delta)^{1/4}(Q(- \Delta)^{1/4} u)-Q(- \Delta)^{1/2}u)\|_{B^0_{1,1}}\\
&&=\sup _{\|h\|_{B^0_{\infty,\infty}}\le 1} \int_{\R^n}\sum_j\sum _{|t-j|\le 3}[(- \Delta)^{1/4}(Q^{j-4}(- \Delta)^{1/4} u_j)-(- \Delta)^{1/2}(Q^{j-4}u_j)h_tdx\nonumber\\
&& =
\sup _{\|h\|_{B^0_{\infty,\infty}}\le 1} \int_{\R^n}\sum_j\sum _{|t-j|\le 3} {\cal{F}}[Q^{j-4}]
{\cal{F}}[(-\Delta)^{1/4}u_j (- \Delta)^{1/4} h_t-(- \Delta)^{1/2}u_jh_t]dx\nonumber\\ && 
 =
\sup _{\|h\|_{B^0_{\infty,\infty}}\le 1} 
\sum_{\ell=0}^\infty\frac{c_{\ell}}{\ell!}\int_{\R^n}
\sum_j\sum _{|t-j|\le 3}
(-i)^{\ell+1}{\cal{F}}[\nabla^{\ell+1}Q^{j-4}]{\cal{F}}[S_{\ell} (- \Delta)^{1/4}u_jh_t)]( x) d x\nonumber\\
&&
\le 
\sup _{\|h\|_{B^0_{\infty,\infty}}\le 1} \|h\|_{B^0_{\infty,\infty}}
\sum_{\ell=0}^\infty\frac{c_{\ell}}{\ell!}\int_{\R^n}
\sum_j  
|\nabla^{\ell+1}(Q)^{j-4}\|S_{\ell} (- \Delta)^{1/4}u_j|d x\nonumber\\
&& 
\le
\sum_{\ell=0}^\infty\frac{c_{\ell}}{\ell!}\int_{\R^n}
(\int_{\R^n}\sum_j |2^{-(k+1/2)j}\nabla^{\ell+1}(Q)^{j-4}\|2^{(\ell+1/2)j}S_{\ell}(- \Delta)^{1/4}u_j|dx\nonumber\\
&& 
\le  C \sum_{\ell=0}^\infty\frac{c_{\ell}}{\ell!}(\int_{\R^n}\sum_j 2^{-2(\ell+1/2)j}|\nabla^{\ell+1} Q^{j-4}|^2 dx)^{1/2}\nonumber\\&&
~~~~~~~~(\int_{\R^n}\sum_j 2^{2(\ell+1/2)j}|S_{\ell} (- \Delta)^{1/4}u_j|^2 dx)^{1/2}\nonumber\\\
&&
\mbox{by Plancherel Theorem}\nonumber \\&&
=
C  \sum_{\ell=0}^\infty\frac{c_{\ell}}{\ell!}
(\int_{\R^n}\sum_j 2^{-2(\ell+1/2)j}|\xi|^{2\ell}|{\cal{F}}[\nabla Q^{j-4}]|^2 d\xi)^{1/2}\nonumber\\
&& ~~~~~~~~(\int_{\R^n}\sum_j 2^{2(\ell+1)j}|\xi|^{-2(\ell+1/2)}|{\cal{F}}[(- \Delta)^{1/4}u_j]|^2 dx)^{1/2}
\nonumber\\\
&&
\le 
C  \sum_{\ell=0}^\infty\frac{c_{\ell}}{\ell!}2^{-3\ell}(\int_{\R^n}\sum_j 2^{-j}|{\cal{F}}[\nabla Q^{j-4}]|^2 d\xi)^{1/2}\nonumber\\
&&~~~~~~~~(\int_{\R^n}\sum_j |{\cal{F}}[(- \Delta)^{1/4}u_j]|^2 dx)^{1/2}
\nonumber\\
&&\le C \sum_{\ell=0}^\infty\frac{c_{\ell}}{\ell!}2^{-3\ell} \|Q\|_{\dot{H}^{1/2}(\R)}\|u\|_{\dot{H}^{1/2}(\R)}\,. ~~~~\Box\nonumber
\end{eqnarray}
The proof of the following Theorems and its localized version can be found in \cite{DL2}\,.
\begin{Theorem}\label{comm2}
Let $u,Q\in \dot{W}^{1/2,q}(\R^n)$ , with $q>2$. Then $T(Q,u), S(Q,u)\in L^{q/2}(\R^n)$  and
\begin{equation}\label{commest4}
\|T(Q,u)\|_{ L^{q/2}}\le C\|(-\Delta)^{1/4}Q\|_{L^q}\|(-\Delta)^{1/4}u\|_{L^q}\,;\end{equation}
\begin{equation}\label{commest5}
\|S(Q,u)\|_{ L^{q/2}}\le C\|(-\Delta)^{1/4}Q\|_{L^q}\|(-\Delta)^{1/4}u\|_{L^q}\,.~\hfill \Box\end{equation}
\end{Theorem}
\begin{Theorem}\label{comm3}
Let $Q\in \dot{H}^{1/2}(\R^n)$ , $u\in \dot{W}^{1/2,q}(\R^n)$  with $q>2$. Then $T(Q,u), S(Q,u)\in L^{\frac{2q}{q+2}}(\R^n)$  and
\begin{equation}\label{commest4}
\|T(Q,u)\|_{ L^{\frac{2q}{q+2}}}\le C\|(-\Delta)^{1/4}Q\|_{L^q}\|(-\Delta)^{1/4}u\|_{L^q}\,;\end{equation}
\begin{equation}\label{commest5}
\|S(Q,u)\|_{ L^{\frac{2q}{q+2}}}\le C\|(-\Delta)^{1/4}Q\|_{L^q}\|(-\Delta)^{1/4}u\|_{L^q}\,.~\hfill \Box\end{equation}
\end{Theorem}
\begin{Remark}
Actually Theorems \ref{comm2} and \ref{comm3} hold for the 2-terms commutators
$$
\tilde T(Q,u)=T(Q,u)-(-\Delta)^{1/4} Q(-\Delta)^{1/4} u=(-\Delta)^{1/4}(Q(-\Delta)^{1/4}u)-Q\Delta^{1/2} u$$
and
$$
\tilde S(Q,u)=S(Q,u)-{\cal{R}}((-\Delta)^{1/4} Q   {\cal{R}}(-\Delta)^{1/4} u)=(-\Delta)^{1/4}[Q(-\Delta)^{1/4} u]-{\cal{R}} (Q\nabla u)\,.~\hfill \Box$$
\end{Remark}
{\bf Aknowledgment.} The author would like to thank Tristan Rivi\`ere for helpful discussions.

 \end{document}